\documentclass[pdflatex,sn-mathphys-num]{sn-jnl}

\usepackage{graphicx}%
\usepackage{subcaption}%
\usepackage{multirow}%
\usepackage{amsmath,amssymb,amsfonts}%
\usepackage{amsthm}%
\usepackage{mathrsfs}%
\usepackage[title]{appendix}%
\usepackage{xcolor}%
\usepackage{textcomp}%
\usepackage{manyfoot}%
\usepackage{booktabs}%
\usepackage{algorithm}%
\usepackage{algorithmicx}%
\usepackage{algpseudocode}%
\usepackage{listings}%
\usepackage{bm}%
\usepackage{onimage}
\usepackage[normalem]{ulem}  
\usepackage{float}
\usepackage{placeins}  

\theoremstyle{thmstyleone}%
\newtheorem{theorem}{Theorem}
\newtheorem{proposition}[theorem]{Proposition}%
\newtheorem{corollary}[theorem]{Corollary}

\theoremstyle{thmstyletwo}%
\newtheorem{remark}{Remark}%

\theoremstyle{thmstylethree}%
\newtheorem{definition}{Definition}%

\newcommand{\ip}[2]{\left\langle #1, #2 \right\rangle}
\newcommand{\U}{\mathcal U}

\begin{document}

\title[Article Title]{Symmetry-reduced model reduction of shift-equivariant systems via operator inference}

\author*[1]{\fnm{Yu} \sur{Shuai}}\email{yu\_shuai@princeton.edu}

\author[1]{\fnm{Clarence W.} \sur{Rowley}}\email{cwrowley@princeton.edu}

\affil*[1]{\orgdiv{Department of Mechanical and Aerospace Engineering}, \orgname{Princeton University}, \orgaddress{\city{Princeton}, \postcode{08540}, \state{New Jersey}, \country{United States}}}

\abstract{We consider data-driven reduced-order models of partial differential equations with shift equivariance.
Shift-equivariant systems typically admit traveling solutions, and the main idea of our approach is to represent the solution in a traveling reference frame, in which it can be described by a relatively small number of basis functions.
Existing methods for operator inference allow one to approximate a reduced-order model directly from data, without knowledge of the full-order dynamics.
Our method adds additional terms to ensure that the reduced-order model not only approximates the spatially frozen profile of the solution, but also estimates the traveling speed as a function of that profile.
We validate our approach using the Kuramoto-Sivashinsky equation, a one-dimensional partial differential equation that exhibits traveling solutions and spatiotemporal chaos.
Results indicate that our method robustly captures traveling solutions, and exhibits improved numerical stability over the standard operator inference approach.
}

\keywords{Shift equivariance, symmetry reduction, reduced-order models, operator inference}


\pacs[MSC Classification (2010)]{35B06, 35Q35, 47A58, 65D15, 68W25}

\maketitle

\section{Introduction}
\label{sec:intro}

Spatially extended systems governed by partial differential equations (PDEs) are central to many areas of physical science, including fluid dynamics, chemical reactions, and atmospheric science.
These PDEs typically describe spatiotemporal scalar, vector, or tensor fields and often exhibit discrete or continuous symmetries such as translation, rotation, or reflection invariance.
In particular, shift-equivariant PDEs, which are invariant under spatial translations, allow solutions to drift continuously in space while preserving their temporal evolution.
This symmetry typically implies that the PDE lacks explicit dependence on spatial coordinates.
In fluid mechanics, for instance, channel flows or circular pipe flows governed
by the Navier-Stokes equations with periodic boundary conditions in the
streamwise and spanwise (azimuthal) directions possess shift equivariance in
both of these directions.
In such systems, propagating coherent structures such as traveling waves, relative periodic orbits and localized turbulent patches have been identified and shown to play central roles in turbulent dynamics~\cite{Wedin2004, Avila2013, Willis2013, Willis2016, Ritter2016, Budanur2017}.
Despite their low-dimensional dynamics, resolving these structures in direct numerical simulations still requires finely resolved spatial grids~\cite{Wedin2004, Avila2013}.
This motivates the development of reduced-order models (ROMs) that exploit the intrinsic low-rank behavior of such systems, providing efficient and interpretable surrogates for full-order models (FOMs).

To construct ROMs for shift-equivariant systems, a common approach is the (Petrov-)Galerkin projection onto spatial modes obtained via proper orthogonal decomposition (POD)~\cite{Sirovich1987, Berkooz1993}.
However, for traveling solutions, the optimal POD modes are Fourier modes~\cite{Berkooz1993, Holmes2012}, which offer limited additional insight, especially when the underlying PDEs are already simulated using Fourier spectral methods.
Moreover, many Fourier modes are typically required to accurately represent propagating solutions, even if the shape of the traveling solution remains relatively constant.
This inefficiency stems from POD’s assumption that spatiotemporal data can be separated into spatial modes and temporal coefficients.
Consequently, the solution manifold has a large Kolmogorov $n$-width, meaning that the projection error onto an $n$-dimensional linear subspace decays slowly~\cite{Greif2019}.

To address this limitation, various methods have been developed to better capture the transport nature of these systems.
Template fitting~\cite{Kirby1992}, also called the method of \emph{slices}~\cite{Rowley2000, Rowley-nl03, Willis2013, Budanur2015}, removes spatial shift by aligning solution snapshots against a ``template'', converting propagating structures into stationary profiles.
These aligned profiles can then be accurately represented with fewer POD modes, and the original solution is recovered via a reconstruction equation that tracks the spatial displacement in real time.
A closely related formulation is the method of \emph{connections}~\cite{Marsden1990, Marsden2000, Rowley2000, Rowley-nl03}, also referred to as \emph{freezing} method~\cite{Beyn2004, Black2020}, where the solution is similarly written as the action of a shift operator on a frozen profile.
A \emph{reconstruction equation} can then be used to reconstruct the translation
amount from the solution in the ``frozen'' frame~\cite{Rowley2000,Beyn2004,Black2020}.
Geometrically, the method of connections restricts the frozen solution to evolve on a state-dependent slice attached to the current frozen profile, whereas the method of slices employs a fixed slice attached to a chosen template~\cite{Rowley2000}.
Both approaches form the foundation of symmetry-reduced, projection-based Galerkin ROMs.

Beyond symmetry reduction, the shifted POD (sPOD) algorithm~\cite{Reiss2019} represents solutions using snapshot-wise shifted POD modes with optimized shifts.
Similarly, the method of unsupervised traveling wave identification with shifting and truncation (UnTWIST)~\cite{Mendible2020} learns the shift as a time-dependent function from a predefined library during training.
While these methods were initially introduced primarily for the offline extraction of low-rank structures, their low-dimensional representations have been successfully extended to construct predictive dynamic ROMs.
For instance, as demonstrated in~\cite{Black2020}, the method of connections can be utilized to dynamically track the continuous shift online, while algorithms like sPOD provide the underlying symmetry-reduced basis.
More generally, alternative approaches circumvent the Kolmogorov $n$-width limitations inherent in fixed global linear projections without relying on an explicit shape-shift decoupling mechanism.
These include ROMs built upon nonlinear manifolds (such as deep convolutional autoencoders~\cite{LeeCarlberg2020}) and time-varying linear projections based on online adaptive-basis framework~\cite{Peherstorfer2020_online_adaptive_basis} where the basis is dynamically updated during time integration via low-rank corrections.

The aforementioned ROMs based on symmetry reduction and online adaptive bases are built upon Galerkin projection.
Consequently, a numerical implementation requires explicit access to the
routines that calculate the governing equations (e.g., the right-hand side of
ODE approximation of the governing PDE); these
are therefore called \emph{intrusive} methods.
To extend model reduction to settings where these routines are inaccessible or treated as a black box, \textit{non-intrusive} data-driven approaches, such as dynamic mode decomposition (DMD)~\cite{Rowley2009, Schmid2010, Tu2012} and operator inference (OpInf)~\cite{Peherstorfer2016, Qian2020, Kramer2024}, offer promising alternatives by relying solely on solution snapshots and their time derivatives.
In the DMD framework, the physics-informed DMD (piDMD) method~\cite{Baddoo2023} enforces commutativity between the learned linear operator and the spatial shift operator to preserve shift equivariance in the reduced model.
However, this constraint restricts the DMD modes to Fourier modes, which often
fail to efficiently capture the coherent structures in the data.
To avoid the reliance on Fourier modes, alternative methods perform dimension reduction on the stationary solution profiles.
The characteristic DMD algorithm~\cite{Sesterhenn2019} eliminates drift by rotating the space-time frame under the assumption of constant group velocity.
More recently, symmetry-reduced DMD (SRDMD)~\cite{Marensi2023, Engel2025} uses template-fitted snapshots as input and constructs a best-fit linear model for the stationary profiles to reveal the low-rank dynamics of coherent structures and to facilitate the identification of unstable periodic orbits embedded in transitional and turbulent flows.
Operator Inference can also be combined with UnTWIST to formulate shifted OpInf ROMs~\cite{Issan2023}, where the shift amount is learned offline as a time-dependent function, and the reduced model is built for the stationary profile using polynomial regression.
Nevertheless, this function can only estimate the shift amounts of unseen solutions via extrapolation rather than dynamically inferring it online from the instantaneous reduced state.
More importantly, these non-intrusive ROMs typically assume polynomial ROM dynamics. 
In contrast, prior intrusive methods such as the symmetry-reduced Galerkin ROM~\cite{Rowley2000} have demonstrated that the dynamics of symmetry-reduced systems contain not only polynomial terms but also a rational symmetry-reduction term which captures the coupling between profile evolution and spatial propagation.
This structure has yet to be fully replicated in non-intrusive ROMs and motivates the framework developed in this work.
Apart from these non-intrusive counterparts of symmetry-reduced Galerkin ROMs, other shift-aware strategies include utilizing deep neural networks built upon a shifted-POD basis (sPOD-NN)~\cite{Burela2025}. 
This method similarly extracts frozen spatial modes but relies on neural networks to map time and parameter inputs to predicted modal amplitudes and shifts. 
Another alternative is the data-driven freezing method via residual minimization~\cite{Cagniart2019}, which determines mode coefficients and shift amounts by minimizing the residual of the full-order dynamics at each discrete timestep. 
While these non-intrusive methods successfully bypass the direct projection of FOM operators, they generally lack a rigorous interpretation of their shift dynamics compared to symmetry reduction techniques. 
Consequently, achieving a non-intrusive reconstruction of the full solution from frozen profiles and shifts remains an open problem.

In this work, we develop a non-intrusive ROM framework for shift-equivariant PDE systems by combining symmetry reduction techniques from~\cite{Rowley2000}, specifically the template fitting process and the reconstruction equation, with Operator Inference (OpInf).
The resulting Symmetry-Reduced Operator Inference (S-R OpInf) ROM serves as a data-driven counterpart to the intrusive S-R Galerkin ROM, yet retains many of its advantages.
Trained on solution data (and possibly the right-hand side dynamical vector field) without access to the individual operator terms in the full-order model (FOM), the S-R OpInf ROM approximates the dynamics of spatially aligned profiles using a compact set of non-Fourier POD modes.
Crucially, our framework enables efficient real-time computation of the shift amount via a reconstruction equation whose coefficients are learned from data.
Moreover, by applying the re-projection technique~\cite{Peherstorfer2020} to pre-process the training snapshots, we show that the S-R OpInf operators can be learned through a strictly convex optimization problem with a unique global optimum matching S-R Galerkin operators.
This establishes a direct correspondence between the intrusive and non-intrusive frameworks and allows insights from the Galerkin formulation to be carried over.
We demonstrate the effectiveness of our method on the Kuramoto–Sivashinsky equation (KSE), focusing on drifting wave solutions with beating temporal dynamics. Specifically, we evaluate the framework across both single- and multiple-trajectory reconstruction tasks.
Our results confirm that, when equipped with suitable penalty regularization, the S-R OpInf ROM accurately reproduces the high-fidelity spatial profiles and shift dynamics.
More importantly, in the multi-solution scenario, our non-intrusive framework demonstrates remarkable robustness on unseen testing data, systematically outperforming its intrusive S-R Galerkin counterpart.
These ingredients also point to potential applications beyond the single-parameter setting studied here---in particular to parameter-dependent and control-oriented reduced-order modeling---which we discuss in Section~\ref{sec:conclusions}.

The remainder of the paper is organized as follows.
Section~\ref{sec:preliminaries} provides a brief overview of shift-equivariant systems and introduces several benchmark ROMs, including the standard Galerkin ROM, OpInf ROM, and the S-R Galerkin ROM.
Section~\ref{sec:sr-opinf} presents the formulation and training algorithm of the proposed S-R OpInf ROM.
In Section~\ref{sec:results}, we apply our method to the Kuramoto–Sivashinsky
equation (KSE) and validate its performance in reconstructing both training and
testing data, and in Section~\ref{sec:conclusions} we summarize our results.

\section{Preliminaries}
\label{sec:preliminaries}

\subsection{Shift-equivariant systems}
\label{sec:shift-equivariant}

We consider a one-dimensional%
\footnote{For clarity of exposition, we present our methodology in the context
  of a one-dimensional shift-equivariant PDE.  This setting captures the
  essential features of the problem while simplifying the notation.  The
  proposed framework, however, generalizes naturally to higher-dimensional
  systems and can be applied to model reduction in those settings without
  conceptual modification.}
partial differential equation (PDE) governing the evolution of a scalar field
$u(x, t)$ defined on the spatial domain $0\leq x \leq L$ and time $t\geq0$ with
periodic boundary conditions $u(x + L, t)\equiv u(x,t)$:

\begin{equation}
    \label{eq:001_governing_pde}
    \partial_t u = f(u),
\end{equation}\\
where $\partial_t u$ denotes the partial derivative of $u(x, t)$ with respect to time $t$, and $f$ is a (spatial) differential operator.
We assume that the initial condition $u(x,0)$ is sufficiently smooth for
existence and uniqueness of solutions to this PDE, and we let $\U$ denote the
space of sufficiently smooth functions of $x$ (i.e., for each $t$,
$u(\cdot,t)\in\U$).

We now define the shift equivariance for system~\eqref{eq:001_governing_pde}:

\begin{definition}[Shift equivariance]
    Let $S_\theta$ denote the shift operator on functions in $\U$:
    \begin{equation}
    \label{eq:002_shift_operator}
    S_{\theta}[v](x) = v(x - \theta),\qquad v\in\U.
    \end{equation}
    The PDE \eqref{eq:001_governing_pde} is shift equivariant if, for all $v\in \U$, $\theta\in\mathbb{R}$, we have

    \begin{equation}
        \label{eq:003_shift_equivariance}
        f(S_\theta[v]) \equiv S_\theta[f(v)].
    \end{equation}
\end{definition}

It follows immediately that, if $f$ is shift equivariant, and $u(x, t)$ is a solution to the PDE~\eqref{eq:001_governing_pde}, then any spatially shifted version $u(x - \theta, t), \forall\theta\in\mathbb{R}$ is also a solution.
In many systems of interest, shift equivariance gives rise to traveling solutions that evolve over time while translating in space.
A traveling solution $u(x, t)$ can be expressed in terms of its spatial profile $u_{\mathrm{fit}}(x, t)$ as follows:
\begin{equation}
    \label{eq:004_traveling_solution}
    u(x, t) = u_{\mathrm{fit}}(x - c(t), t),
\end{equation}\\
where $c(t)$ is the shift amount.
For instance, for the simple advection equation
\begin{equation*}
  \partial_t u + a\, \partial_x u = 0,
\end{equation*}
the corresponding $f(u)=-a\partial_x u$ is shift equivariant, and the equation admits traveling wave solutions of the form $u(x,t)=w(x-c(t))$, with $c(t)=at$.
More generally, shift-equivariant systems can also support traveling solutions with time-varying profiles and shift speeds.

\subsection{Standard reduced-order models}
\label{sec:roms}

We now give an overview of two standard model reduction techniques: the Galerkin
method, and the Operator Inference method.  Both approaches seek to
approximate the full-state vector ${u}(t)\in\U$ using a low-dimensional basis,
but they differ in how the reduced dynamics are obtained.  In the Galerkin
method, the dynamics are obtained through projection of the individual components in the right-hand side $f$ of~\eqref{eq:001_governing_pde}, whereas in Operator
Inference, the model is learned from data via regression.

Regarding notation, for the remainder of this paper, we regard the
PDE~\eqref{eq:001_governing_pde} as a dynamical system evolving in the space~$\U$, which
consists of (periodic) functions of space~$x$.  Thus, at a particular time~$t$,
$u(t)\in\U$ is a function of space, equivalent to $u(\cdot,t)$.

\subsubsection{Galerkin models for quadratic dynamics}
\label{sec:galerkin}

To construct a Galerkin model, one usually begins by collecting an ensemble of
snapshots $\{{u}(t_m)\}_{m = 0}^{N_t-1}$ (each a function of space), which serves
as the training dataset.  One then chooses $n$ basis functions
$\varphi_1,\ldots,\varphi_n$ for a subspace of
$\U$ on which to project the equations.  We will assume that $\U$ is an inner
product space, with the standard inner product
\begin{equation}
  \label{eq:005_inner_product}
  \ip{v}{w} = \frac{1}{L} \int_0^L v(x)w(x)\,dx,\qquad v,w\in\U,
\end{equation}
with induced norm $\|v\| = \sqrt{\ip{v}{v}}$, and that the basis functions are orthonormal (i.e., $\ip{\varphi_i}{\varphi_j}=\delta_{ij}$).
We then approximate $u(t)$ by the reduced expansion
\begin{equation}
    \label{eq:006_reduced_expansion}
    {u}(t) \approx u_n(t) = \sum_{i = 1}^{n}a_i(t){\varphi}_i,
\end{equation}\\
where the coefficients ${a_i}(t)\in\mathbb{R}$ are the reduced variables.

A common choice for the basis is to perform proper orthogonal decomposition
(POD)~\cite{Sirovich1987} on the set of snapshots
$\{u(t_0),\ldots, u(t_{N_t-1})\}$.  This is the procedure we shall
use in the examples in Section~\ref{sec:results}, but the developments below are
independent of the choice of basis; we assume only that the $\{\varphi_i\}$
are orthonormal.

For most of what follows, we will assume that the differential operator $f$
in~\eqref{eq:001_governing_pde} is \emph{quadratic} with no constant term: that is, it has the
form
\begin{equation}
  \label{eq:007_quadratic_fom}
  f(u) = A u + B(u,u),
\end{equation}
where $A:\U\to\U$ is linear, and $B:\U\times\U\to\U$ is bilinear
(linear in each argument) and symmetric (i.e., $B(u,v)=B(v,u)$ for all
$u,v\in\U$).  A nonzero constant term, if present, can be eliminated by a change
of variables, so we omit it without loss of generality.

Galerkin projection then proceeds as follows.  We would like to choose the
coefficients $a_i(t)$ such that the expansion~\eqref{eq:006_reduced_expansion}
satisfies the original PDE~\eqref{eq:001_governing_pde}.  However, this
is in general not possible, as the true dynamics might drive $u(t)$ out of the subspace
spanned by our $n$ basis functions.  So, in Galerkin projection, we simply project
the dynamics (orthogonally) onto this subspace, and choose $a_i(t)$ so that
\begin{equation}
  \label{eq:008_galerkin_projection}
  \ip{\partial_t u_n}{\varphi_i} = \ip{f(u_n)}{\strut\varphi_i},\qquad i=1,\ldots,n.
\end{equation}
Using the fact that the basis functions $\varphi_i$ are orthonormal, we obtain
\begin{equation*}
  \dot a_i = \ip{Au_n + B(u_n,u_n)}{\strut\varphi_i},
\end{equation*}
and inserting the expression~\eqref{eq:006_reduced_expansion} for $u_n$, we obtain
the Galerkin model
\begin{equation}
  \label{eq:009_galerkin_model}
  \dot a_i = A_{ij} a_j + B_{ijk} a_j a_k,\qquad i=1,\ldots,n,
\end{equation}
with
\begin{equation}
  \label{eq:010_galerkin_coeffs}
  \begin{aligned}
    A_{ij}  & = \ip{A\varphi_j}{\strut\varphi_i} \\
    B_{ijk} & = \ip{B(\varphi_j,\varphi_k)}{\strut\varphi_i},
  \end{aligned}
\end{equation}
where in~\eqref{eq:009_galerkin_model} we use the convention that repeated indices are summed.
Observe that, because of the symmetry of $B$, it follows that $B_{ijk}=B_{ikj}$.
Also note that, in order to calculate the coefficients $A_{ij}$ and
$B_{ijk}$, one needs to be able to evaluate the linear and bilinear components of $f$, and thus this
method is considered ``intrusive.''

\subsubsection{Operator Inference models}
\label{sec:opinf}

The formulation of Galerkin models described in the previous section requires
explicit access to the operators of the full-order dynamics~\eqref{eq:007_quadratic_fom}, in
order to calculate the coefficients in~\eqref{eq:010_galerkin_coeffs}.  However, in many practical
scenarios, the full model may be given as a ``black box'' providing only the
function $f$ (not its individual components $A$ and $B$), or simply a
timestepper that advances $u(t)$ to $u(t+\Delta t$).

To build reduced-order models without requiring access to the full-order
operators, the operator inference method~\cite{Peherstorfer2016} offers
a natural and flexible approach.  The initial steps are similar to those of the
Galerkin method, as we determine a reduced basis ${\varphi}_i$ from snapshots
and project the full-state trajectories onto this basis to obtain reduced
variables ${a_i}(t_m) = \ip{u(t_m)}{\varphi_i}$.  We again assume the
dynamics are quadratic, as in~\eqref{eq:007_quadratic_fom}, but without explicit access to these
operators, we simply postulate the quadratic form of the reduced-order model as
in~\eqref{eq:009_galerkin_model}, and determine the coefficients $A_{ij}$ and $B_{ijk}$
by fitting to the available data.

In order to determine these coefficients, we will also require the time
derivatives of the snapshots, $\partial_t u(t_m)$.  If we have access to the
function $f$, these may be computed by simply evaluating $f(u(t_m))$.  If not,
and if the snapshots are spaced closely enough in time, we may approximate the
time derivatives using finite differences:
\begin{equation*}
  f\big(u(t_m)\big) \approx \frac{u(t_{m+1}) - u(t_m)}{t_{m+1}-t_m}.
\end{equation*}
From our dataset $\{u(t_m)\}$, we then define
\begin{equation}
  \label{eq:011_opinf_data}
  a_i(t_m) = \ip{u(t_m)}{\strut\varphi_i},\qquad
  f_i(t_m) = \ip{f\big(u(t_m)\big)}{\varphi_i},
\end{equation}
and the model parameters $\mathcal{O}=(A_{ij},B_{ijk})$ are determined from the minimization problem
\begin{equation}
  \label{eq:012_opinf_lstsq}
  \min_{\mathcal{O}} \sum_{m = 0}^{N_t-1}\sum_{i=1}^n
  \Big[A_{ij} a_j(t_m) + B_{ijk} a_j(t_m) a_k(t_m) - f_i(t_m)\Big]^2 + R(\mathcal{O}),
\end{equation}
where $R(\mathcal{O})$ denotes a regularization function that penalizes large
values of the coefficients.
Given that the parameters in~$\mathcal{O}$ appear
linearly inside the square brackets, if the regularization function is also
quadratic in the parameters, then this is a least-squares problem that is
\emph{strictly} convex---and thus has a unique global minimizer---provided the regression matrix has full column rank, for instance with enough columns or with Tikhonov regularization~\cite{Peherstorfer2016}.

\begin{remark}[Enforcing index symmetry in the data-driven quadratic operator]
\label{rem:opinf-idx-sym}
To make the inferred quadratic operator consistent with the Galerkin tensor in structure, we
also enforce the index symmetry $B_{ijk}=B_{ikj}$.
Thus only $n^2(n+1)/2$ independent quadratic coefficients are actually learned, rather than $n^3$.
In practice, we include each quadratic monomial $a_j a_k$ only once for $1\le j\le k\le n$ during the OpInf regression.
After the regression, the fitted coefficient of $a_j^2$ is
assigned to $B_{ijj}$, while the fitted coefficient of $a_j a_k$ with $j<k$ is
divided equally between $B_{ijk}$ and $B_{ikj}$. 
Thus, the compact expression $B_{ijk}a_j a_k$ should be understood as the
duplicate-free quadratic expansion in the offline learning stage, and as the equivalent
symmetrized full-tensor contraction during online evaluation.
Unless otherwise stated, all quadratic coefficient tensors written below in compact tensor form follow this duplicate-free/symmetrized convention.
\end{remark}

\subsection{Symmetry reduction and the Galerkin method}
\label{sec:sr-galerkin}

The performance of reduced-order models, whether formulated via
Galerkin projection or regression, depends heavily on the choice of basis
functions~${\varphi}_i$.
It is well known that for systems with periodic or translational symmetry, the optimal (in the sense of the $L^2$ projection error) POD basis consists of Fourier modes~\cite{Berkooz1993}.  In order to represent
localized features, typically a large number of Fourier modes is required, so
reduced-order models with a small number of Fourier modes usually offer limited
insight into the underlying coherent structures.

To address this limitation, a natural approach is to first remove the spatial
translation from the evolution of solutions and then compute POD modes from the
resulting spatially aligned solution profiles~\cite{Kirby1992}.  This strategy
allows for a compact representation with fewer modes, from which reduced-order
models can be derived to capture the essential dynamics of the spatial profiles.
Once an approximation to the profile is obtained, the full solution can be
reconstructed by reintroducing the spatial shift determined from the model.
This methodology is the essence of the symmetry reduction technique~\cite{Kirby1992, Rowley2000}, also known as the method of slices~\cite{Rowley-nl03, Willis2013, Budanur2015}.  The main contribution of this
paper is to incorporate symmetry reduction into the operator inference method,
and this will described in Section~\ref{sec:sr-opinf}.  Before addressing this,
we first give an overview of symmetry reduction with the Galerkin method.

\subsubsection{Reducing the shift-equivariance symmetry}
\label{sec:reducing-symmetry}

Symmetry-reduced models are based on representing a solution $u(t)$ as a shifted
version of a function $u_{\mathrm{fit}}(t)$, and this amounts to writing the equations of
motion in a moving reference frame in which the translation has been removed.
In terms of the shift operator~$S_c$ defined in
Section~\ref{sec:shift-equivariant}, we write
\begin{equation}
  \label{eq:013_shifted_representation}
  u(t) = S_{c(t)} u_{\mathrm{fit}}(t).
\end{equation}
Assume for the moment that we have determined an appropriate shift amount $c(t)$
(we will discuss the choice of $c(t)$ in detail below).  Then the dynamics of
$u_{\mathrm{fit}}(t)$ must satisfy
\begin{align*}
  \partial_t u_{\mathrm{fit}}(x,t) &= \frac{d}{dt}\big[u(x+c(t),t)\big]\\
  &= \partial_t u(x+c(t),t) + \dot c(t)\,\partial_x u(x+c(t),t),
\end{align*}
where $\dot c$ denotes $dc/dt$.  Using the dynamics $\partial_t u = f(u)$ and writing this in terms of the shift
operator, we obtain
\begin{align*}
  \partial_t u_{\mathrm{fit}}(t) &= S_{-c(t)} f(u(t)) + \dot c(t) S_{-c(t)}\partial_x u(t)\\
  &= f(S_{-c(t)} u(t)) + \dot c(t)\, \partial_x S_{-c(t)} u(t),
\end{align*}
where we have used shift equivariance of $f$.  Thanks to~\eqref{eq:013_shifted_representation}, this then
becomes
\begin{equation}
  \label{eq:014_comoving_dynamics}
  \partial_t u_{\mathrm{fit}} = f(u_{\mathrm{fit}}) + \dot c\, \partial_xu_{\mathrm{fit}}.
\end{equation}
These are the dynamics in the moving reference frame; observe that the
equations~\eqref{eq:014_comoving_dynamics} depend only on $\dot c$, not on $c$
itself.

We now return to the question of how to determine the shift amount $c(t)$.  At
each time~$t$, we will choose the shift amount $c(t)$ so that the shifted
solution $u_{\mathrm{fit}}(t)$ is optimally ``aligned'' with a pre-chosen ``template
function'' $u_0\in\U$.  That is, we will choose
\begin{equation}
  \label{eq:015_template_fitting}
  c(t) = \arg\min_c \|S_{-c} u(t) - u_0\|^2,
\end{equation}
where the norm on $\U$ is the standard one induced by the inner
product~\eqref{eq:005_inner_product}.  This is equivalent to
\begin{align*}
  c(t) & = \arg\max_c \ip{S_{-c} u(t)}{u_0} \\
       & = \arg\max_c \ip{u(t)}{S_c u_0}.
\end{align*}
A necessary condition for a local maximum is therefore
\begin{equation*}
  \frac{d}{dc}\ip{u(t)}{S_c u_0} = 0.
\end{equation*}
A quick calculation reveals that, for any $v\in\U$,
\begin{equation*}
  \frac{d}{dc} S_c[v](x) = \frac{d}{dc} v(x-c) = -\partial_x v(x-c) = -S_c[\partial_x v](x),
\end{equation*}
and thus the condition for a local maximum becomes
\begin{align*}
  \ip{u(t)}{-S_{c(t)} \partial_x u_0} &= 0\\
  \iff \ip{S_{-c(t)} u(t)}{\partial_x u_0} &= 0\\
  \iff \ip{u_{\mathrm{fit}}(t)}{\partial_x u_0} &= 0.
\end{align*}
As explained in~\cite{Rowley2000}, this condition has an intuitive geometric
interpretation: the shifted solution $u_{\mathrm{fit}}$ must lie in a subspace of~$\U$ orthogonal
to the function $\partial_x u_0$ (this subspace is called a \emph{slice}, hence
the \emph{method of slices} mentioned earlier).

Since this relation must hold for $u_{\mathrm{fit}}(t)$ at every time~$t$, we may
differentiate with respect to~$t$ to obtain
\begin{equation*}
  \frac{d}{dt} \ip{\strut u_{\mathrm{fit}}(t)}{\partial_x u_0} = 0,
\end{equation*}
and substituting the dynamics~\eqref{eq:014_comoving_dynamics}, we have
\begin{equation*}
  \ip{\strut f(u_{\mathrm{fit}}) + \dot c\, \partial_x u_{\mathrm{fit}}}{\partial_x u_0} = 0.
\end{equation*}
Thus, solving for $\dot c$, we obtain
\begin{equation}
  \label{eq:016_shift_reconstruction}
  \dot c = -\frac{\ip{\strut f(u_{\mathrm{fit}})}{\partial_x u_0}}{\ip{\strut\partial_x u_{\mathrm{fit}}}{\partial_x u_0}}.
\end{equation}
This expression then determines $\dot c$ in the dynamics~\eqref{eq:014_comoving_dynamics}.
Furthermore, once a solution $u_{\mathrm{fit}}(t)$ to~\eqref{eq:014_comoving_dynamics} is obtained, the shift
amount $c(t)$ can be determined by integrating~\eqref{eq:016_shift_reconstruction} in time, and used
to reconstruct the solution~$u(t)$ in the original reference
frame.  We therefore call~\eqref{eq:016_shift_reconstruction} a ``reconstruction equation,''
after~\cite{Rowley2000}.

The shift reconstruction equation used here should be distinguished from those in the method of connections or the freezing method~\cite{Beyn2004, Cagniart2019, Black2020}.
In these template-free frameworks, the shift speed is determined via minimizing the norm of the time derivative of $u_{\mathrm{fit}}(t)$, which gives $\dot c = -\ip{\strut f(u_{\mathrm{fit}})}{\partial_x u_{\mathrm{fit}}}/\ip{\strut\partial_x u_{\mathrm{fit}}}{\partial_x u_{\mathrm{fit}}}$~\cite{Beyn2004}, or the residual error of the time derivative of $u_{\mathrm{fit}}(t)$ after substituting the approximation ~\eqref{eq:017_comoving_expansion} introduced in the next section~\cite{Cagniart2019, Black2020}.
While these methods remove the instantaneous spatial-translation component in the dynamics of the current state,
the frozen profile may still exhibit a slow residual drift due to the accumulated geometric phase~\cite{Rowley-nl03}.
In addition, replacing $\partial_x u_0$ in~\eqref{eq:016_shift_reconstruction} with $\partial_x u_{\mathrm{fit}}$ also raises the polynomial order of the numerator and denominator of the shift dynamics by 1, which increases both the online evaluation cost
and the regression cost for a non-intrusive symmetry-reduced ROM to be developed.
We therefore adopt the method of slices using a fixed external template to derive our reconstruction equation.

\subsubsection{Symmetry-reduced Galerkin for quadratic dynamics}
\label{sec:sr-galerkin-quadratic}

We now combine the symmetry-reduction approach with the Galerkin
method from Section~\ref{sec:galerkin}.  We write $u(t)=S_{c(t)}u_{\mathrm{fit}}(t)$, with
\begin{equation}
  \label{eq:017_comoving_expansion}
  u_{\mathrm{fit}}(t) \approx u_{\mathrm{fit},n}(t) = \sum_{i=1}^n a_i(t)\varphi_i,
\end{equation}
where the orthonormal basis functions $\varphi_i\in\U$ are computed from the template-fitted snapshots $\{u_{\mathrm{fit}}(t_m)\}$ instead of the raw snapshots $\{u(t_m)\}$.
\footnote{It is worth noting that we also explored the conventional mean-subtracted expansion $u_{\mathrm{fit}}(t) \approx u_{\mathrm{fit},n}(t) = \bar{u} + \sum_{i=1}^n a_i(t)\varphi_i$, where $\bar{u}$ denotes the time average of $u_{\mathrm{fit}}(t)$. However, ROMs constructed upon this formulation yielded poor performance. Specifically, we found that relying on a fixed $\bar{u}$ pre-computed from the training data introduces a static bias, which degrades the reconstruction of testing trajectories with different time-averaged profiles.
To achieve greater generalizability across multiple distinct trajectories, we ultimately opted to omit the mean subtraction.}
The time evolution of $u_{\mathrm{fit}}(t)$ is then
determined by~\eqref{eq:014_comoving_dynamics}, and
Galerkin projection gives
\begin{equation}
  \label{eq:018_sr_galerkin_projection}
  \ip{\partial_tu_{\mathrm{fit},n}}{\strut\varphi_i} =
  \ip{f(u_{\mathrm{fit},n}) + \dot c\,\partial_xu_{\mathrm{fit},n}}{\strut\varphi_i},
\end{equation}
where $\dot c$ is given by~\eqref{eq:016_shift_reconstruction}.  Inserting our
expansion~\eqref{eq:017_comoving_expansion}, and using orthogonality of $\varphi_i$, we obtain
\begin{equation}
  \label{eq:019_sr_galerkin_adot}
  \dot a_i = A_{ij} a_j + B_{ijk} a_j a_k + \dot c\, C_{ij} a_j,
\end{equation}
for $i=1,\ldots,n$, where we again employ the summation convention, and
coefficients are given by~\eqref{eq:010_galerkin_coeffs}, together with
\begin{equation}
  \label{eq:020_recon_coeffs}
  \begin{aligned}
    C_{ij} & = \ip{\partial_x\varphi_j}{\strut\varphi_i},
  \end{aligned}
\end{equation}
and $\dot c$ has the form
\begin{equation}
  \label{eq:021_sr_galerkin_cdot}
  \dot c = -\frac{p_j a_j + Q_{jk} a_j a_k}{s_j a_j},
\end{equation}
with
\begin{equation}
  \label{eq:022_shift_coeffs}
  \begin{aligned}
    p_j    & = \ip{\strut A\varphi_j}{\partial_x u_0} \\
    Q_{jk} & = \ip{\strut B(\varphi_j,\varphi_k)}{\partial_x u_0}             \\
    s_j    & = \ip{\strut\partial_x\varphi_j}{\partial_x u_0}
  \end{aligned}
\end{equation}
Observe that there is a symmetry in the coefficients $Q_{jk}=Q_{kj}$, because of
the symmetry of~$B$.  Like the Galerkin method discussed in
Section~\ref{sec:galerkin}, this method is intrusive.
Our objective in the next section is to obtain a non-intrusive
method.

\section{Symmetry-reduced operator inference}
\label{sec:sr-opinf}

We now describe the main contribution of this paper, a method combining symmetry
reduction with operator inference.  We first describe a baseline approach with a
na\"{\i}ve non-intrusive approximation of the shift speed $\dot{c}(t)$, and show
its limitations.  We then describe our proposed method, which incorporates
additional operators into the reduced-order model.

Here, the term ``non-intrusiveness'' is used in the operator-inference sense: the reduced operators are learned without explicit evaluation or projection of the full-order operators, their adjoints, or their Jacobians.
However, our framework still performs basic grid-based operations on the full-state data vectors, such as shifting with a given shift operator, computing inner products and spatial derivatives, and state resetting (see Appendix~\Ref{app:sr}).
Nevertheless, these are strictly data-processing steps.
Provided that the spatial grid is sufficiently resolved, these routines depend solely on the available full-order state vectors and remain independent of the underlying physical equations and spatial discretization details (e.g., spectral-based versus interpolation-based shifting of snapshots). Thus, they preserve the non-intrusive, data-driven nature of the ROM framework.

\subsection{Basic formulation and its drawbacks}
\label{sec:sr-opinf-basic}

In standard operator inference (Section~\ref{sec:opinf}), one learns
the coefficients $A_{ij}$ and $B_{ijk}$ of reduced dynamics, from which we assemble
\begin{equation}
  \label{eq:023_reduced_rhs}
  \dot a_i = A_{ij} a_j + B_{ijk} a_j a_k,
\end{equation}
so that the dynamics of the reduced state shares the same algebraic structure as the intrusive Galerkin counterpart.
Combining this ansatz with the symmetry reduction described in Section~\ref{sec:sr-galerkin} gives a symmetry-reduced operator inference ROM with the same form as the symmetry-reduced Galerkin ROM in~\eqref{eq:019_sr_galerkin_adot}:
\begin{equation*}
  \dot a_i = A_{ij} a_j + B_{ijk} a_j a_k + \dot c\, C_{ij} a_j.
\end{equation*}
Here, the matrix $C$ is computed as in~\eqref{eq:020_recon_coeffs}, and the tensors $A$ and $B$ are coefficients to be learned.
To construct $C$ and infer $A$ and $B$, we work throughout with
the template-fitted snapshots.
The POD basis $\{\varphi_i\}_{i=1}^n$ is computed
from the shifted profiles $\{u_{\mathrm{fit}}(t_m)\}_{m=0}^{N_t-1}$, and we project the shifted snapshots and their velocities onto this
basis to obtain the reduced data $a_i(t_m) = \ip{u_{\mathrm{fit}}(t_m)}{\varphi_i}$ and
$f_i(t_m) = \ip{f(u_{\mathrm{fit}}(t_m))}{\varphi_i}$.  The coefficients $A_{ij}$, $B_{ijk}$
are then learned from these data exactly as in the standard operator-inference
least-squares problem~\eqref{eq:012_opinf_lstsq}.

It remains to determine the dynamics of the shift amount $c$ in the ROM.
Following Section~\ref{sec:sr-galerkin}, given the reduced state $a=(a_1,\ldots,a_n)$, we reconstruct the shifted profile and its velocity
in the full space and substitute them into the reconstruction
equation~\eqref{eq:016_shift_reconstruction}.
The profile is reconstructed as
$u_{\mathrm{fit},n}(a)=\sum_i a_i\varphi_i$.
For the velocity $f(u)$, a direct evaluation $f(u_{\mathrm{fit},n}(a))$ would require evaluating the full-order right-hand side online,
making the ROM expensive.
To circumvent this, we replace it by its
orthogonal projection onto the basis: $f(u_{\mathrm{fit},n}(a))\approx \sum_i\ip{f(u_{\mathrm{fit},n}(a))}{\varphi_i}\varphi_i$.
By the formulation of symmetry-reduced operator inference, we can evaluate this inner product using the previously learned coefficients: $\ip{f(u_{\mathrm{fit},n}(a))}{\varphi_i}\approx A_{ij} a_j + B_{ijk} a_j a_k$.
Therefore, we define
\begin{equation}
  \label{eq:024_naive_reconstruction}
  f_{\mathrm{fit},n}(a) = \sum_{i=1}^n\Bigg(\sum_{j=1}^n A_{ij} a_j +\sum_{j, k=1}^n B_{ijk} a_j a_k\Bigg)\varphi_i,\qquad
  u_{\mathrm{fit},n}(a)= \sum_{i=1}^n a_i\varphi_i,
\end{equation}
and take the shift speed to be computed from these projection-based approximations
\begin{equation}
  \label{eq:025_naive_shift_speed}
  \dot c = -\frac{\ip{f_{\mathrm{fit},n}(a)}{\partial_x u_0}}
  {\ip{\partial_x u_{\mathrm{fit},n}(a)}{\partial_x u_0}}.
\end{equation}
The numerator is quadratic in $a$ and inherits the learned coefficients $A_{ij}$,
$B_{ijk}$; the denominator is linear in $a$, with known coefficients given
by~\eqref{eq:022_shift_coeffs}.
The right-hand side of~\eqref{eq:025_naive_shift_speed} is thus evaluated as a fixed-coefficient rational function of the ROM state $a$.

This procedure seems reasonable, but it has a subtle yet serious flaw: the functions $f(u_{\mathrm{fit}})$ often do
not lie in the subspace that $u_{\mathrm{fit}}$ lies in, and this can result in
drastically incorrect approximations of the shift speed~$\dot c$.
As an illustrative example, consider the advection-diffusion equation
\begin{equation}
  \label{eq:026_advection_diffusion}
  \partial_t u + v\,\partial_x u = \partial_x^2 u,
\end{equation}
in a periodic domain.
The accurate shift speed determined by
equation~\eqref{eq:016_shift_reconstruction} is
\begin{equation*}
  \dot c = -\frac{\ip{-v\,\partial_x u_{\mathrm{fit}} + \partial_x^2 u_{\mathrm{fit}}}{\partial_x
      u_0}}{\ip{\partial_x u_{\mathrm{fit}}}{\partial_x u_0}}
  = v + \frac{\ip{\partial_x^2u_{\mathrm{fit}}}{\partial_x u_0}}{\ip{\partial_x u_{\mathrm{fit}}}{\partial_x u_0}}.
\end{equation*}
Now, suppose the initial condition is an even function $u(x,0)=g(x)$, where
$g(-x)=g(x)$, and take the template to be the initial condition itself, $u_0=u(x,0)$.
Since the initial snapshot coincides with the template,
template fitting~\eqref{eq:015_template_fitting} returns $c(0)=0$ exactly, and $u_{\mathrm{fit}}(x,0)=u(x,0)$ is even.
Then
$\partial_x u_0$ is odd and $\partial_x^2 u_{\mathrm{fit}}$ is even, so the numerator in
the last term
above is zero, and the ``true'' shift speed is simply $\dot c = v$.
Furthermore, the dynamics of $u_{\mathrm{fit}}$ as given by~\eqref{eq:014_comoving_dynamics} are
\begin{equation}
  \label{eq:027_diffusion_reduced}
  \partial_t u_{\mathrm{fit}} = \partial_x^2 u_{\mathrm{fit}},
\end{equation}
a simple diffusion equation.  The shifting procedure has removed the advection part of the
equation, as desired.

Now, compare this with the approximation~\eqref{eq:025_naive_shift_speed}.  Since the initial condition is an even function, the solution~$u_{\mathrm{fit}}(t)$
of~\eqref{eq:027_diffusion_reduced} will remain even for all time (as $\partial_x^2u_{\mathrm{fit}}$ is
even), so our basis functions~$\varphi_i$ (for instance, determined by POD) will
also be even functions.
Thus, if we determine $f_{\mathrm{fit},n}$ according
to~\eqref{eq:024_naive_reconstruction}, it will be even as well.
But then in our approximation of $\dot c$ in~\eqref{eq:025_naive_shift_speed}, the numerator consists of the inner product between an even function ($f_{\mathrm{fit}, n}$) and an odd function
($\partial_x u_0$).  Since even and odd functions are orthogonal, the shift
speed predicted by our model will always be zero!

One solution could be to consider a larger subspace, which captures $f(u_{\mathrm{fit}})$
as well as~$u_{\mathrm{fit}}$: for instance, if one uses POD to determine the basis
functions $\varphi_i$, one could include snapshots of $f(u(t_m))$ in addition to
snapshots of $u(t_m)$.  However, this would result in a larger-dimensional
subspace, and a less ``efficient'' reduced-order model.  Below we take a
different approach.

\subsection{An improved version of Symmetry-Reduced Operator Inference}
\label{sec:sr-opinf-improved}

As we have seen, approximating the shift speed~$\dot c$ according
to~\eqref{eq:025_naive_shift_speed} has serious drawbacks.  To address these, we propose learning a
model for $\dot c$ at the same time that we learn the model for $f_i$
in~\eqref{eq:023_reduced_rhs}.
In particular, since we know the true $\dot c$ has the form~\eqref{eq:021_sr_galerkin_cdot}, we
consider a model of the same form
\begin{equation}
  \label{eq:028_sr_opinf_shift_model}
  \dot c(a) = -\frac{p_j a_j + Q_{jk} a_j a_k}{s_j a_j},
\end{equation}
in which the coefficients $s_j$ are computed from~\eqref{eq:022_shift_coeffs},
while the coefficients $p_j$ and $Q_{jk}$ need to be learned from the data.
Because of the symmetry $Q_{jk}=Q_{kj}$, the number of additional parameters to
learn is $n + n(n+1)/2$.

For a given dataset $\{u(t_m)\}$, we begin by aligning the data to a template
function $u_0\in\U$, to obtain shifted snapshots $u_{\mathrm{fit}}(t_m)$ and shift amounts $c(t_m)$ for our dataset.  As with standard operator
inference~\eqref{eq:011_opinf_data}, we project the (shifted) snapshots and their time
derivatives, to define
\begin{equation}
  \label{eq:029_sr_opinf_data}
  a_i(t_m) = \ip{u_{\mathrm{fit}}(t_m)}{\strut\varphi_i},
  \qquad f_i(t_m) = \ip{f\big(u_{\mathrm{fit}}(t_m)\big)}{\varphi_i}.
\end{equation}
Due to the shift equivariance of $f$, the following relation holds:
\begin{equation}
  \label{eq:030_fhat_equivariance}
  f\big(u_{\mathrm{fit}}(t_m)\big) = f\big(S_{-c(t_m)} u(t_m)\big) = S_{-c(t_m)}f\big(u(t_m)\big),
\end{equation}
which ensures that as long as the derivative data of the original solution $f(u(t_m))$ are available, we can evaluate $f\big(u_{\mathrm{fit}}(t_m)\big)$ without requiring access to the function $f$.
Furthermore, by defining
\begin{equation}
  \label{eq:031_sr_opinf_rdef}
  g(t_m) = \ip{f\big(u_{\mathrm{fit}}(t_m)\big)}{\partial_x u_0},
\end{equation}
we can then formulate two separate minimization problems for $(A_{ij},B_{ijk})$ and $(p_j,Q_{jk})$, respectively:
\begin{subequations}
  \label{eq:032_sr_opinf_loss}
  \begin{align}
    \min_{A,B}\ &\sum_{m=0}^{N_t-1}\sum_{i=1}^n\bigg[A_{ij}a_j(t_m) + B_{ijk}a_j(t_m)a_k(t_m) - f_i(t_m)\bigg]^2 + \lambda_{\mathrm{poly}}\big(\|A\|_F^2 + \|B\|_F^2\big), \label{eq:032a_poly_loss}\\
    \min_{p,Q}\ &\sum_{m=0}^{N_t-1}\bigg[p_j a_j(t_m) + Q_{jk}a_j(t_m)a_k(t_m) - g(t_m)\bigg]^2 + \lambda_{\mathrm{num}}\big(\|p\|_F^2 + \|Q\|_F^2\big). \label{eq:032b_num_loss}
  \end{align}
\end{subequations}
Here $\|\cdot\|_F$ denotes the Frobenius norm, which is the square root of the sum of the squares of all elements in the matrix or tensor.
The minimization problem~\eqref{eq:032a_poly_loss} is the standard operator-inference
least-squares problem~\eqref{eq:012_opinf_lstsq} for the polynomial
coefficients with penalty weight $\lambda_{\mathrm{poly}}\geq 0$, and~\eqref{eq:032b_num_loss} fits the numerator
coefficients $(p_j,Q_{jk})$ of the reconstruction
equation~\eqref{eq:028_sr_opinf_shift_model} with another penalty $\lambda_{\mathrm{num}}\geq 0$.
Each subproblem is a Tikhonov-regularized least-squares problem.
A nonzero penalty ensures that the objective function is strictly convex, thereby guaranteeing a unique minimizer even when the data matrix is rank-deficient.
While the approach does involve more optimization parameters than the
original operator inference problem, the number of additional parameters ($n + n(n+1)/2$) is small compared with the number of the parameters in the original
operator inference problems (which is more than $n^2(n+1)/2$).

The final reduced-order model is then given by
\begin{equation}
  \label{eq:033_sr_opinf_rom}
  \dot a_i = A_{ij} a_j + B_{ijk} a_ja_k + \dot c(a)\, C_{ij}a_j,
\end{equation}
with $\dot c(a)$ given by~\eqref{eq:028_sr_opinf_shift_model} and coefficients learned from~\eqref{eq:032_sr_opinf_loss}.
Algorithm~\ref{alg:sr-opinf} summarizes all of the steps involved in our
symmetry-reduced operator inference procedure, which we call S-R OpInf.

\begin{algorithm}
  \caption{The learning algorithm of Symmetry-Reduced Operator Inference}\label{alg:sr-opinf}
  \begin{algorithmic}[1]
    \Require Snapshots of the solution ${u}(t_m)\in\U$ and the velocities
    $f(u(t_m))\in\U$, $m=1,\ldots,N_t$. The template function $u_0\in\mathcal{U}$ for template fitting. The dimension of the desired reduced-order model $n$.
    \Ensure The Symmetry-Reduced Operator Inference model~\eqref{eq:033_sr_opinf_rom} and~\eqref{eq:028_sr_opinf_shift_model} describing the time evolution of the reduced state ${a}(t)=(a_1(t),\ldots,a_n(t))$ and the shift amount $c(t)$.
    \For{$m = 1, \dots, N_t$}
      \State Perform template fitting to obtain the shifted snapshots $u_{\mathrm{fit}}(t_m)$ and the corresponding velocities $f(u_{\mathrm{fit}}(t_m))$ according to~\eqref{eq:013_shifted_representation}, ~\eqref{eq:015_template_fitting}, and ~\eqref{eq:030_fhat_equivariance}.
    \EndFor
    \State Perform POD on the set of snapshots $\{u_{\mathrm{fit}}(t_m)\}$,
    to obtain orthonormal basis functions $\varphi_i\in\U$, $i=1,\ldots,n$.
    \For{$m = 1, \ldots, N_t$}
      \State Determine $a_i(t_m)$ and $f_i(t_m)$ from~\eqref{eq:029_sr_opinf_data}, and
      $g(t_m)$ from~\eqref{eq:031_sr_opinf_rdef}.
    \EndFor
    \State Compute the coefficients $s_j$ from~\eqref{eq:022_shift_coeffs}
    and $C_{ij}$ from~\eqref{eq:020_recon_coeffs}.
    \State Learn the coefficients $A_{ij}$, $B_{ijk}$, $p_j$, and
    $Q_{jk}$ by solving the optimization problem~\eqref{eq:032_sr_opinf_loss} (noting the
    symmetry, $B_{ijk}=B_{ikj}$ and $Q_{jk}=Q_{kj}$).
    \State The reduced-order model is given by~\eqref{eq:033_sr_opinf_rom} with $\dot c$ given by~\eqref{eq:028_sr_opinf_shift_model}.
  \end{algorithmic}
\end{algorithm}

\section{Application}
\label{sec:results}
In this work, we apply the proposed S-R OpInf ROM to perform model reduction on a representative shift-equivariant PDE system: the Kuramoto–Sivashinsky equation (KSE), given by 
\begin{equation}
    \label{eq:034_kse}
    \partial_t u = f(u) = -u\,\partial_x u - \partial_x^2 u - \nu\,\partial_x^4 u,
\end{equation}
on the domain $x\in[0, 2\pi]$ with periodic boundary conditions.
The function $f$ is thus of the form~\eqref{eq:007_quadratic_fom} with
\begin{equation*}
Au = -\partial_x^2 u - \nu \partial_x^4 u,\qquad B(u,v) = -\frac{1}{2}(u\,\partial_x v + v\,\partial_x u).
\end{equation*}

The linear part of the KSE consists of a destabilizing second-order term $-\partial_x^2 u$ and a stabilizing fourth-order hyper-diffusion term $-\nu\,\partial_x^4 u$.
This combination implies that low-wavenumber modes are linearly unstable, while
high-wavenumber modes are linearly damped. As a result, the KSE exhibits
dynamics qualitatively similar to the energy cascade in fluid turbulence governed by the Navier–Stokes equations, and is often referred to as a model for “one-dimensional turbulence”~\cite{Holmes2012}.
Another key reason for the widespread use of the KSE in reduced-order modeling lies in its rich solution dynamics, especially in the parameter regime $\nu\in(4/89, 2/43)$ where the system exhibits a variety of solutions with low-order dynamics~\cite{Kevrekidis1990}.
In this work, we focus on the globally shifting, beating-wave periodic solutions as a particular class of such solutions and collect their snapshots as the training data for our ROMs.
These solutions make the KSE an ideal testbed for evaluating ROMs designed to handle high-dimensional shift-equivariant systems with spatial drifting solutions.

We apply S-R OpInf to two online-reconstruction tasks of increasing difficulty.  The first is to reconstruct a single base traveling-beating solution, shown in Fig.~\ref{fig:original_and_template_fitted_training_snapshots}.  The second is to reconstruct a whole family of traveling-beating solutions whose initial transients are generated by linearly unstable perturbations of the base solution, which probes how well the framework generalizes across trajectories.  Section~\ref{sec:numerics} details the numerical setup common to both.

\subsection{Numerical details}
\label{sec:numerics}
We simulate~\eqref{eq:034_kse} at $\nu = 4/87$ using a Fourier pseudospectral method on a $256$-point grid applying Orszag's $3/2$ rule to dealias the quadratic nonlinearity~\cite{Orszag1971}.
We advance the solution with a semi-implicit timestepper (RK3CN,~\cite{Zhang1985, Peyret2002}), using Crank--Nicolson for the linear terms and a third-order explicit Runge--Kutta scheme for the nonlinear term, at a step $\Delta t_{\mathrm{FOM}} = 10^{-3}$.

To obtain beating-wave solutions, we first simulate the FOM from the harmonic initial condition
\begin{equation}
    \label{eq:035_kse_harmonics_as_base_ic_starter}
    u(x,0) = -\sin x + 2\cos 2x + 3\cos 3x - 4\sin 4x
\end{equation}
for $120$ time units to settle onto the beating-wave attractor~\cite{Kirby1992, Rowley2000}, and denote the resulting state by $u_{\mathrm b}(x) = u(x, t = 120)$.
We then build the datasets for the two of  our numerical tasks with $u_b(x)$.
For the task of reconstructing a single base solution, we simulate the FOM with $u_{\mathrm b}(x)$ as the new initial condition $u(x, 0) = u_{\mathrm b}(x)$ over the next $10$ time units, collecting $N_t = 101$ temporally equispaced snapshots at a sample interval $\Delta t_{s} = 0.1$.
These snapshots form the dataset $\{u(t_m), f(u(t_m))\}_{m = 0}^{N_t-1}$ (the dependence on $x$ is omitted for brevity) with $t_0 = 0$, where we assume direct access to the function $f$ as a whole in this work.
For the numerical experiment on multiple-trajectory reconstruction, we generate an ensemble of initial conditions $u(x,0) = u_{\mathrm b}(x) + u_{\mathrm p}(x)$.
The perturbation field $u_{\mathrm p}(x)$ is constructed as
\begin{equation*}
  u_{\mathrm p}(x) = \gamma\,\|u_{\mathrm b}\|\,\frac{p(x)}{\|p\|},\qquad
  p(x)=\sum_{k=2}^{k_{\max}}\big(a_k\cos kx + b_k\sin kx\big),
\end{equation*}
where $\gamma$ denotes the relative perturbation intensity.
For each trajectory realization, the coefficients $a_k$ and $b_k$ are drawn independently from a standard normal distribution (i.i.d.\ $\mathcal{N}(0,1)$). 
We restrict the summation up to $k_{\max}=\lfloor 1/\sqrt{\nu}\,\rfloor = 4$, which corresponds to the largest linearly unstable wavenumber $k$ with a positive linear growth rate $k^2(1-\nu k^2)$.
Specifically, we draw 10 training trajectories at $\gamma =10\%$ and generate 30 independent trajectories for each intensity level $\gamma \in \{12, 13, 14, 15, 16\}\%$ as testing set.
Every trajectory is simulated for 10 time units from a distinct initial condition, yielding 101 equispaced snapshots per realization.

\begin{figure}[tbp]
  \centering
  \begin{subfigure}[b]{0.32\textwidth}
    \centering
    \begin{tikzonimage}[width=\linewidth]{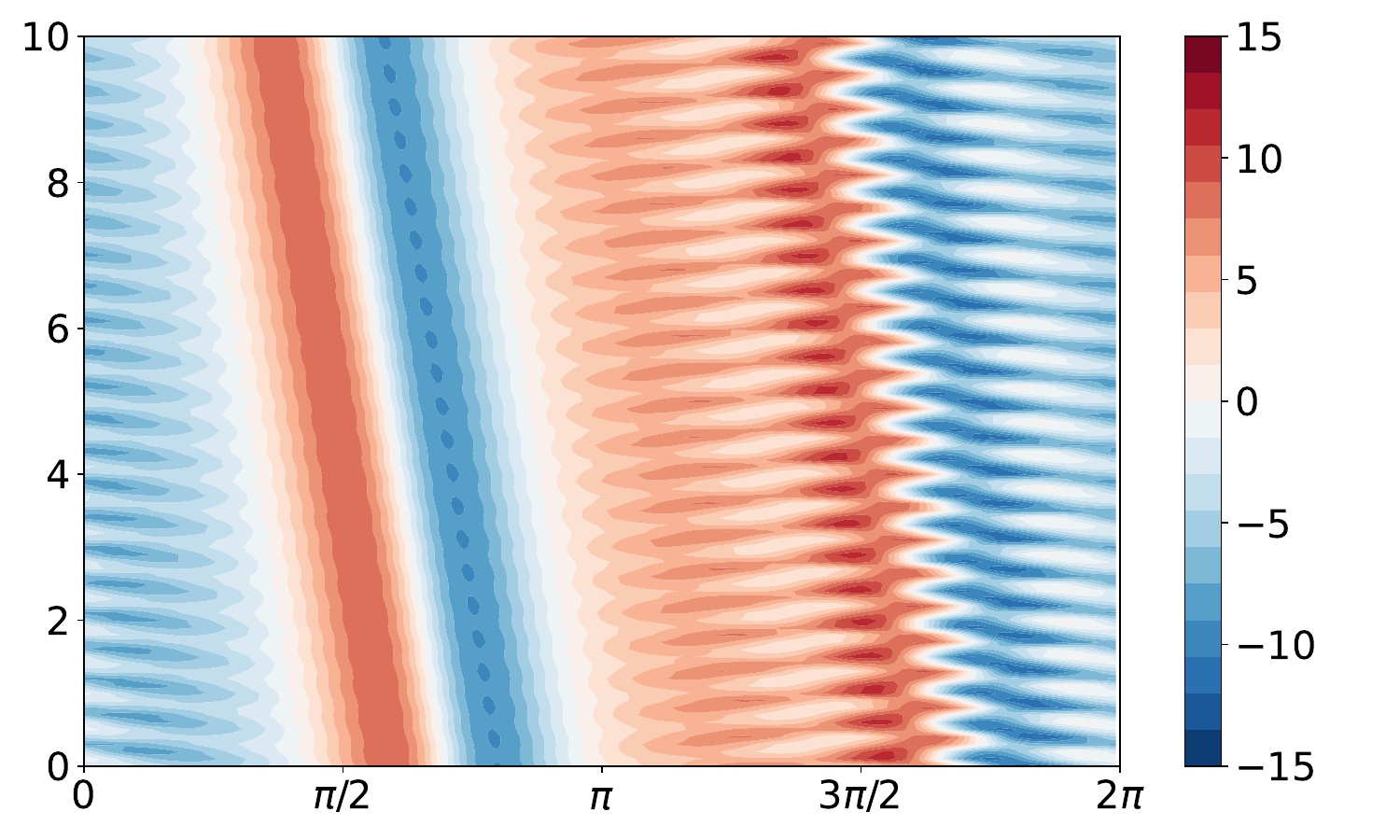}
      \node at (0.45, -0.05) {$x$};
      \node[overlay] at (-0.08, 0.52) {$t$};
      \node at (0.45, -0.25) {(a)};
    \end{tikzonimage}
    \captionsetup{labelformat=empty}
    \caption{\mbox{}}
    \vspace{-1.5em}
    \label{fig:original_and_template_fitted_training_snapshots_panel_a}
  \end{subfigure}
  \hfill
  \begin{subfigure}[b]{0.32\textwidth}
    \centering
    \begin{tikzonimage}[width=\linewidth]{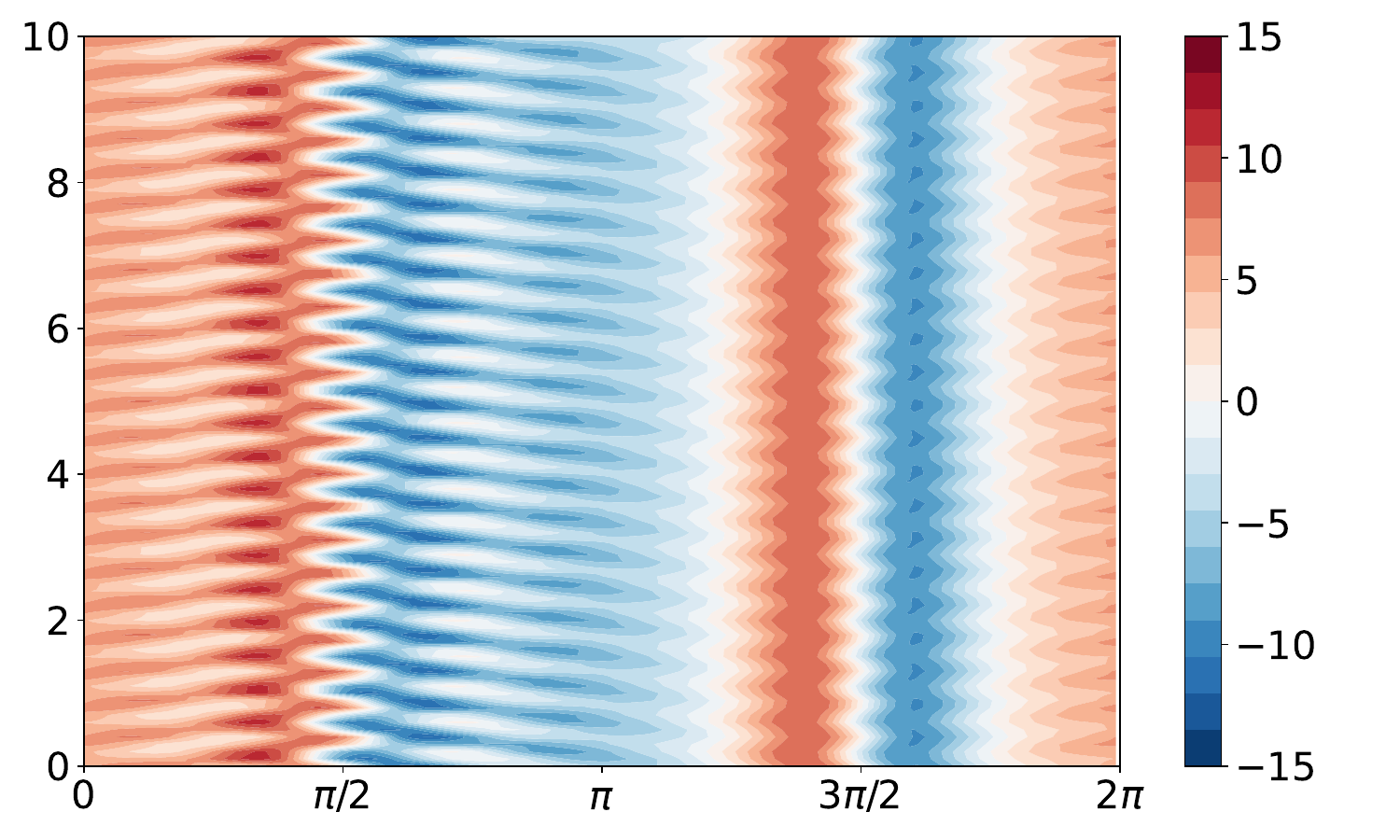}
      \node at (0.45, -0.05) {$x$};
      \node at (0.45, -0.25) {(b)};
    \end{tikzonimage}
    \captionsetup{labelformat=empty}
    \caption{\mbox{}}
    \vspace{-1.5em}
    \label{fig:original_and_template_fitted_training_snapshots_panel_b}
  \end{subfigure}
  \hfill
  \begin{subfigure}[b]{0.32\textwidth}
    \centering
    \begin{tikzonimage}[width=\linewidth]{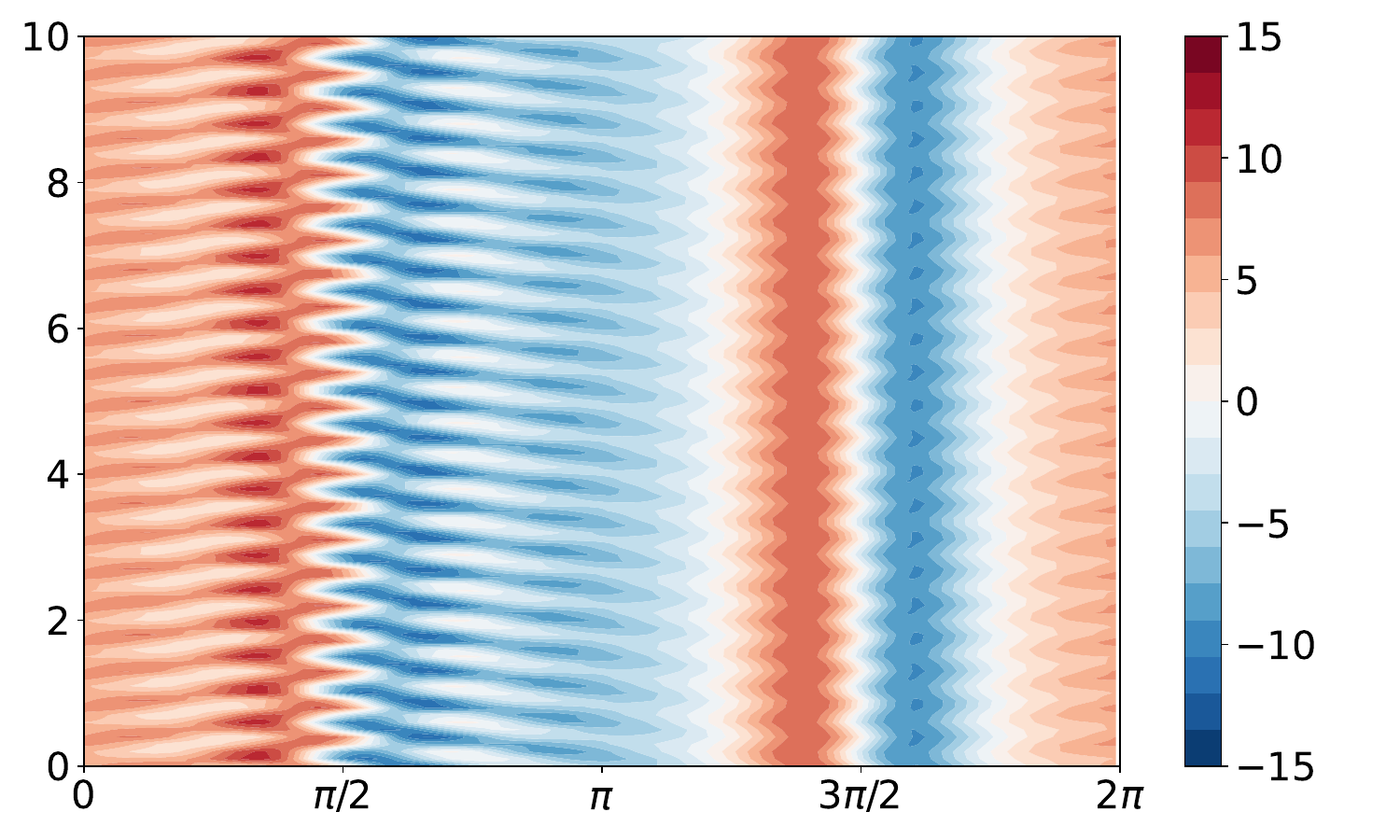}
      \node at (0.45, -0.05) {$x$};
      \node at (0.45, -0.25) {(c)};
    \end{tikzonimage}
    \captionsetup{labelformat=empty}
    \caption{\mbox{}}
    \vspace{-1.5em}
    \label{fig:original_and_template_fitted_training_snapshots_panel_c}
  \end{subfigure}
  \caption{The base full-order solution at $\nu = 4/87$ on a 256-point grid, exhibiting an upstream-traveling beating pattern.  Panel~(a) shows $u(t)$ in the original frame; panel~(b) shows the template-fitted profile $u_{\mathrm{fit}}(t)$ with the drift removed; and panel~(c) shows the template fitting carried out with a cubic-spline-interpolated shift operator.  Color denotes the value of $u$ on the spatial grid $[0,2\pi]$.}
  \label{fig:original_and_template_fitted_training_snapshots}
\end{figure}

We align the snapshots to the template $u_0(x) = \cos x$, the real part of the first Fourier mode~\cite{Budanur2015}.
On the periodic domain $[0,2\pi]$, the template-fitting problem~\eqref{eq:015_template_fitting} has the closed-form maximizer $c(t) = \arg\!\big(\langle u(t),\cos x\rangle + \mathrm{i}\,\langle u(t),\sin x\rangle\big)$ with ``$\mathrm{i}$'' the imaginary unit, from which we obtain a continuous shift (unique modulo $2\pi$~\cite{Marensi2023}) and form the template-fitted snapshots $u_{\mathrm{fit}}(t_m) = S_{-c(t_m)}u(t_m)$ and $f(u_{\mathrm{fit}}(t_m))$.
The shift operator $S_{-c}$ is applied exactly in Fourier space by multiplying the $k$-th coefficient by $\exp(\mathrm{i}kc)$ following the definition~\eqref{eq:002_shift_operator}.
On a physical grid, however, the shift, spatial derivatives, and inner products are approximated by interpolation, differentiation, and quadrature.
This approximation introduces spatial discretization errors, including an $O(\Delta x^p)$ contribution from a $p$-th order interpolation-based shift.
Our supplementary tests, in which the template-fitted snapshots and the
resulting symmetry-reduced ROMs are computed on 40-point and 256-point physical
grids using cubic-spline interpolation, show that the discrepancy from exact spectral
shifting is visible at 40 points but negligible at 256 points.
Fig.~\ref{fig:original_and_template_fitted_training_snapshots}(b)--(c)
compares the fitted snapshots for the 256-point case, showing no evident difference.
We therefore recommend using
sufficiently fine spatial resolution when applying this non-intrusive model reduction
method to problems for which only discrete measurements are available.

We build the reduced basis by POD of the template-fitted snapshots.
For the single-solution task, $n = 5$ modes capture $99.99\%$ of the snapshot energy, while for the multiple-solution case we retain $n = 20$ modes.
In both experiments, once the POD basis is fixed, we compute the reduced dataset $\{a_i(t_m), f_i(t_m), g(t_m)\}_{m = 0}^{N_t-1}$ from~\eqref{eq:029_sr_opinf_data} and~\eqref{eq:031_sr_opinf_rdef} for the S-R OpInf least-squares problems in~\eqref{eq:032_sr_opinf_loss}.
The Tikhonov penalty weights are selected by grid search, giving $(\lambda_{\mathrm{poly}}, \lambda_{\mathrm{num}}) = (1.56\times10^{-1}, 6.90\times10^{-6})$ for the single-solution case and $(\lambda_{\mathrm{poly}}, \lambda_{\mathrm{num}}) = (10^{-13}, 10^{-7})$ for the multiple-solution case.
Each of the two least-squares problems~\eqref{eq:032a_poly_loss}--\eqref{eq:032b_num_loss} is small and is solved directly, without iteration.

Once trained, each ROM is evaluated by integrating the reduced state $a(t)$ and the shift $c(t)$ simultaneously, using $y(t) = [a_i(t); c(t)]\in\mathbb R^{n+1}$ as the augmented state governed by the autonomous system
\begin{equation*}
  \frac{d}{dt}\begin{bmatrix} a_i\\ c\end{bmatrix}
  = \begin{bmatrix} A_{ij} a_j + B_{ijk} a_j a_k + \dot c(a)\,C_{ij} a_j\\[2pt] \dot c(a)\end{bmatrix},
  \qquad \dot c(a) = -\big(p_j a_j + Q_{jk} a_j a_k\big)\frac{s_j a_j}{(s_j a_j)^2 + \delta}.
\end{equation*}
Here, $A_{ij}$, $B_{ijk}$, $p_j$, $Q_{jk}$ are learned from~\eqref{eq:032_sr_opinf_loss}, while $C_{ij}$, and $s_j$ are computed from~\eqref{eq:020_recon_coeffs} and~\eqref{eq:022_shift_coeffs}.
To initialize the time integration, the initial shift $c(0)$ is determined by template fitting the full-state initial condition $u(0)$ via~\eqref{eq:015_template_fitting}.
Subsequently, the initial ROM state $a(0)$ is obtained by projecting this aligned initial condition $u_{\mathrm{fit}}(0) = S_{-c(0)}u(0)$ onto the low-dimensional subspace $a_i(0) = \ip{u_{\mathrm{fit}}(0)}{\varphi_i}$.
Note that during online evaluation, we use a safeguarded form of the reconstruction equation for $\dot{c}(a)$ instead of its original form~\eqref{eq:028_sr_opinf_shift_model} by replacing the factor $1/(s_ja_j)$ with $(s_ja_j)/\big((s_ja_j)^2 + \delta\big)$. Here, $\delta\geq 0$ is a hyperparameter introduced to prevent numerical blow-up caused by division by zero, which could otherwise occur if $s_ja_j$ becomes anomalously small due to ROM approximation errors.
In all numerical experiments reported here, however, the denominator $s_j a_j$ stays well away from zero, so we simply set $\delta = 0$.
It is also worth mentioning that the rational term in $\dot c$ can make the system stiff.
Therefore, we integrate the equations using the implicit Radau~IIA method~\cite{Hairer1996} with an adaptive timestep (\texttt{scipy.integrate.solve\_ivp}), and set the relative and absolute tolerances to $10^{-3}$ and $10^{-6}$, respectively.
If the adaptive step falls below $10^{-5}$ or the state diverges, the ROM will be flagged as numerically unstable and the simulation will be terminated.

To evaluate the performance of the S-R OpInf ROMs, we introduce several benchmark ROMs for the single-solution and the multiple-solution experiments for comparison.
In the first task we compare four ROMs in total: the intrusive S-R Galerkin ROM; the S-R OpInf ROM trained with zero penalties; the S-R OpInf ROM trained with the grid-searched optimal penalties; and the S-R OpInf ROM trained on the re-projected dataset with zero penalties.
In particular, the \emph{re-projected} dataset consists of preprocessed training snapshots, which was introduced for OpInf of discrete-time systems in~\cite{Peherstorfer2020}.
Under suitable conditions it guarantees that the learned S-R OpInf reduced-order operators replicate exactly those obtained from Galerkin projection.
In this work we extend this guarantee to our symmetry-reduced, continuous-time setting and summarize the construction of the re-projected dataset in Appendix~\ref{app:reprojection}.
For the multiple-solution case we compare two ROMs: the S-R Galerkin ROM and the S-R OpInf ROM trained with the grid-searched optimal penalties.
In neither case do we include the standard, non-symmetry-reduced Galerkin or OpInf ROMs, which are known to perform poorly here---an $8$-mode Galerkin ROM is significantly outperformed by a $3$-mode S-R Galerkin ROM~\cite{Rowley2000}.

\subsection{Reconstructions of the single traveling solution}
\label{sec:single-reconstruction}

\begin{figure}[tbp]
  \centering
  \begin{subfigure}[b]{0.46\textwidth}
    \centering
    \begin{tikzonimage}[width=\linewidth]{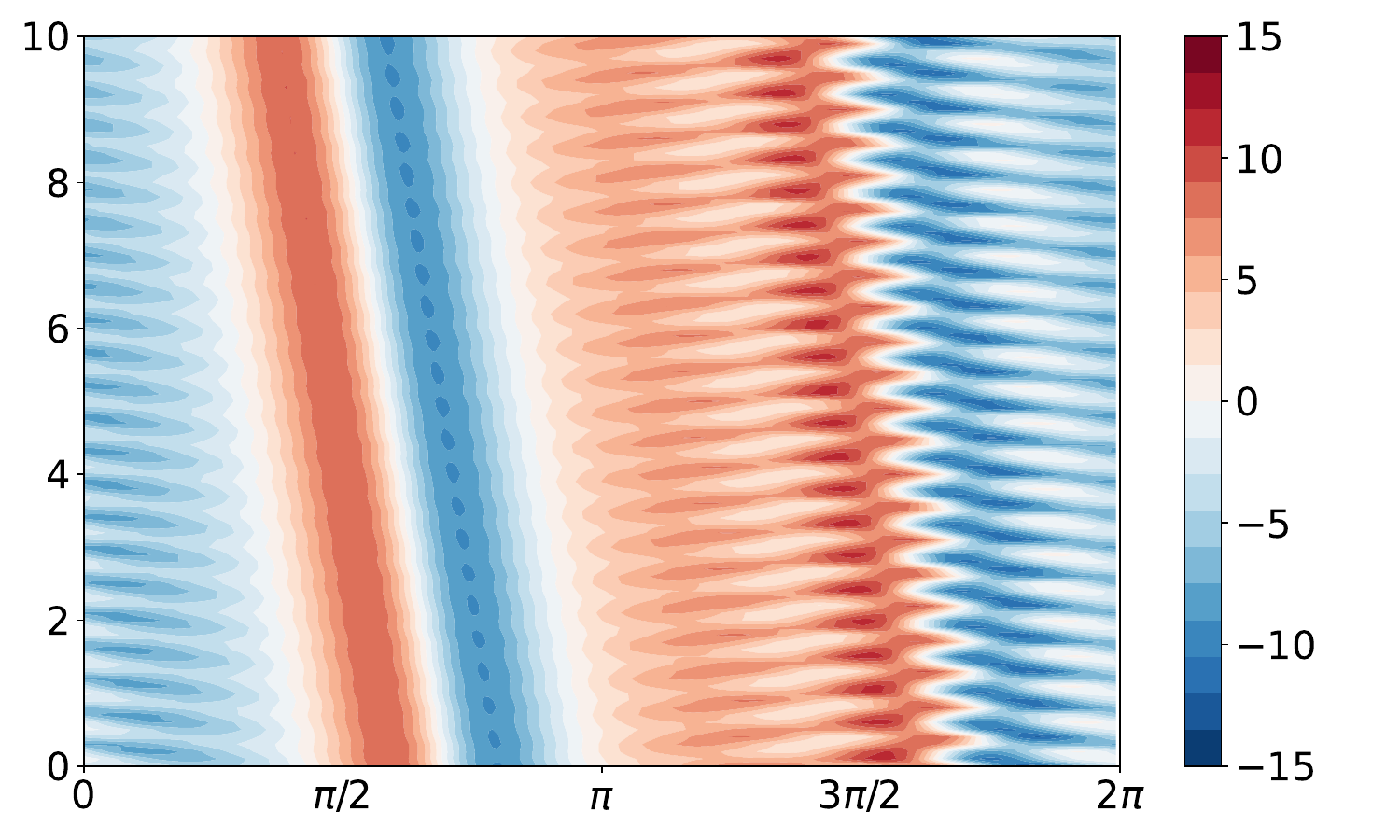}
      \node[anchor=south] at (0.45, 1.0) {\large S-R Galerkin, 6.12\%};
      \node at (0.45, -0.05) {$x$};
      \node[overlay] at (-0.08, 0.52) {$t$};
      \node at (0.45, -0.18) {(a)};
    \end{tikzonimage}
    \captionsetup{labelformat=empty}
    \caption{\mbox{}}
    \vspace{-1.5em}
    \label{fig:ROM_training_snapshots_panel_a}
  \end{subfigure}
  \hfill
  \begin{subfigure}[b]{0.46\textwidth}
    \centering
    \begin{tikzonimage}[width=\linewidth]{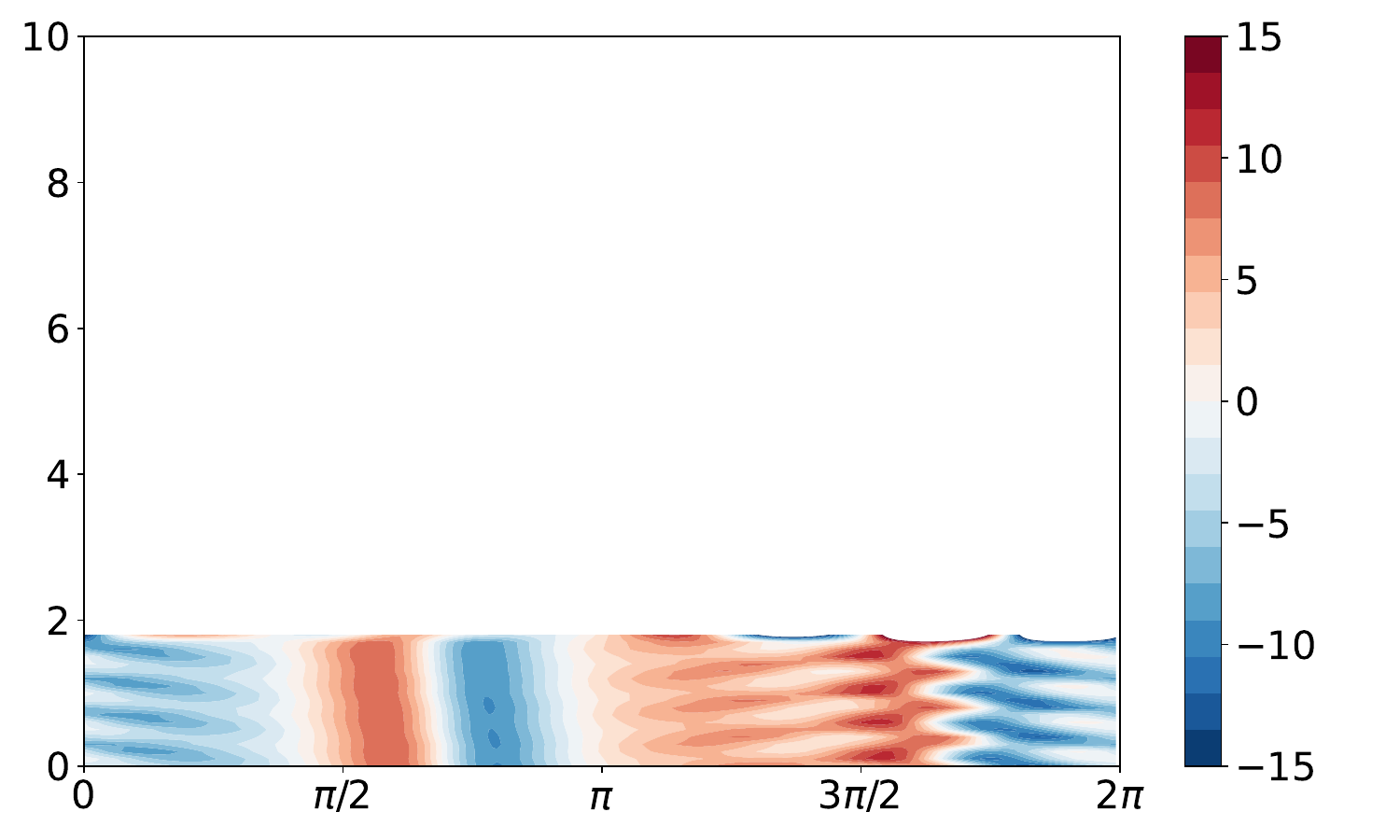}
      \node[anchor=south] at (0.45, 1.0) {\large S-R OpInf w/o reproj};
      \node at (0.45, -0.05) {$x$};
      \node at (0.45, 0.6) {(blow-up)};
      \node at (0.45, -0.18) {(b)};
    \end{tikzonimage}
    \captionsetup{labelformat=empty}
    \caption{\mbox{}}
    \vspace{-1.5em}
    \label{fig:ROM_training_snapshots_panel_b}
  \end{subfigure}
  \begin{subfigure}[b]{0.46\textwidth}
    \centering
    \begin{tikzonimage}[width=\linewidth]{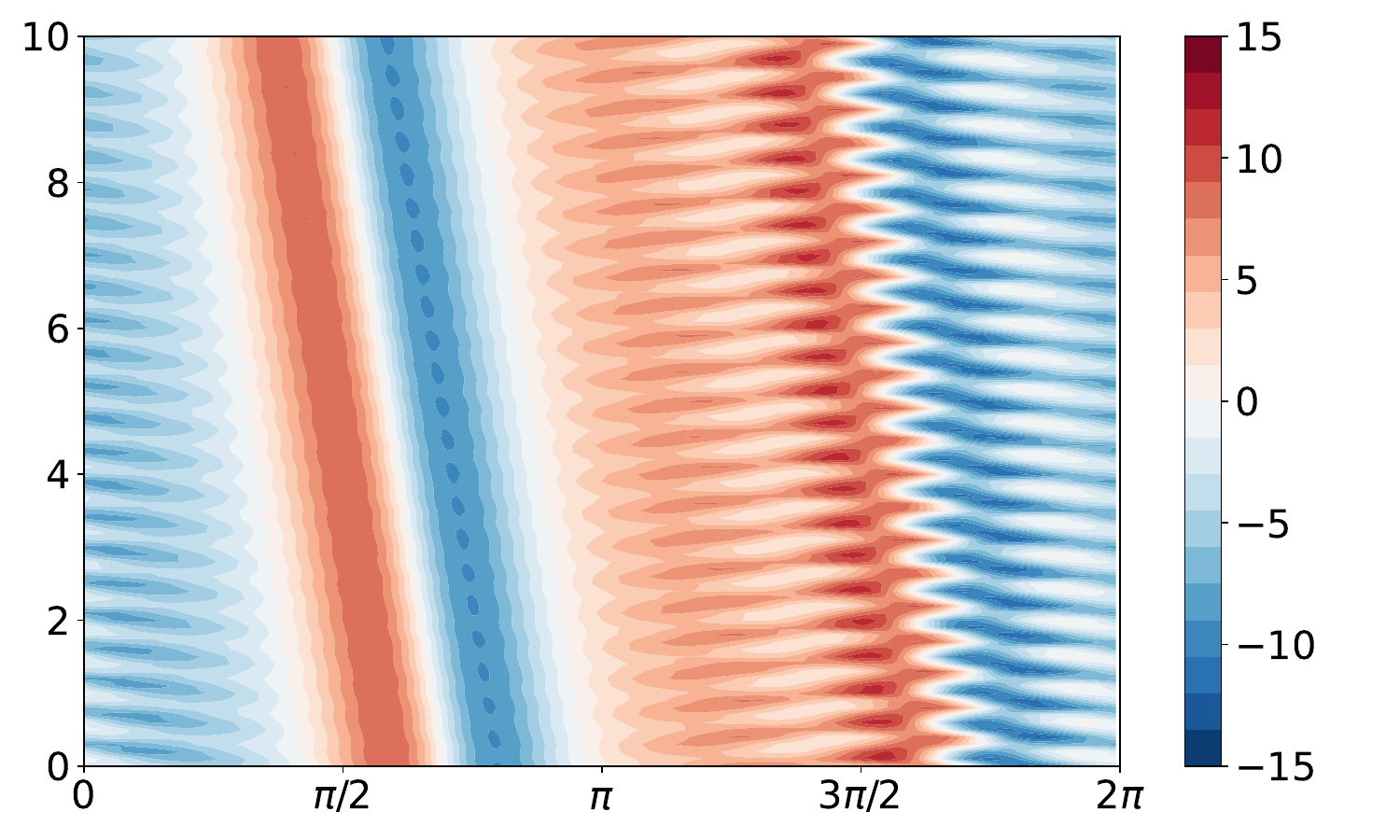}
      \node[anchor=south] at (0.45, 1.0) {\large S-R OpInf w/ reproj, 6.12\%};
      \node at (0.45, -0.05) {$x$};
      \node[overlay] at (-0.08, 0.52) {$t$};
      \node at (0.45, -0.18) {(c)};
    \end{tikzonimage}
    \captionsetup{labelformat=empty}
    \caption{\mbox{}}
    \vspace{-1.5em}
    \label{fig:ROM_training_snapshots_panel_c}
  \end{subfigure}
  \hfill
  \begin{subfigure}[b]{0.46\textwidth}
    \centering
    \begin{tikzonimage}[width=\linewidth]{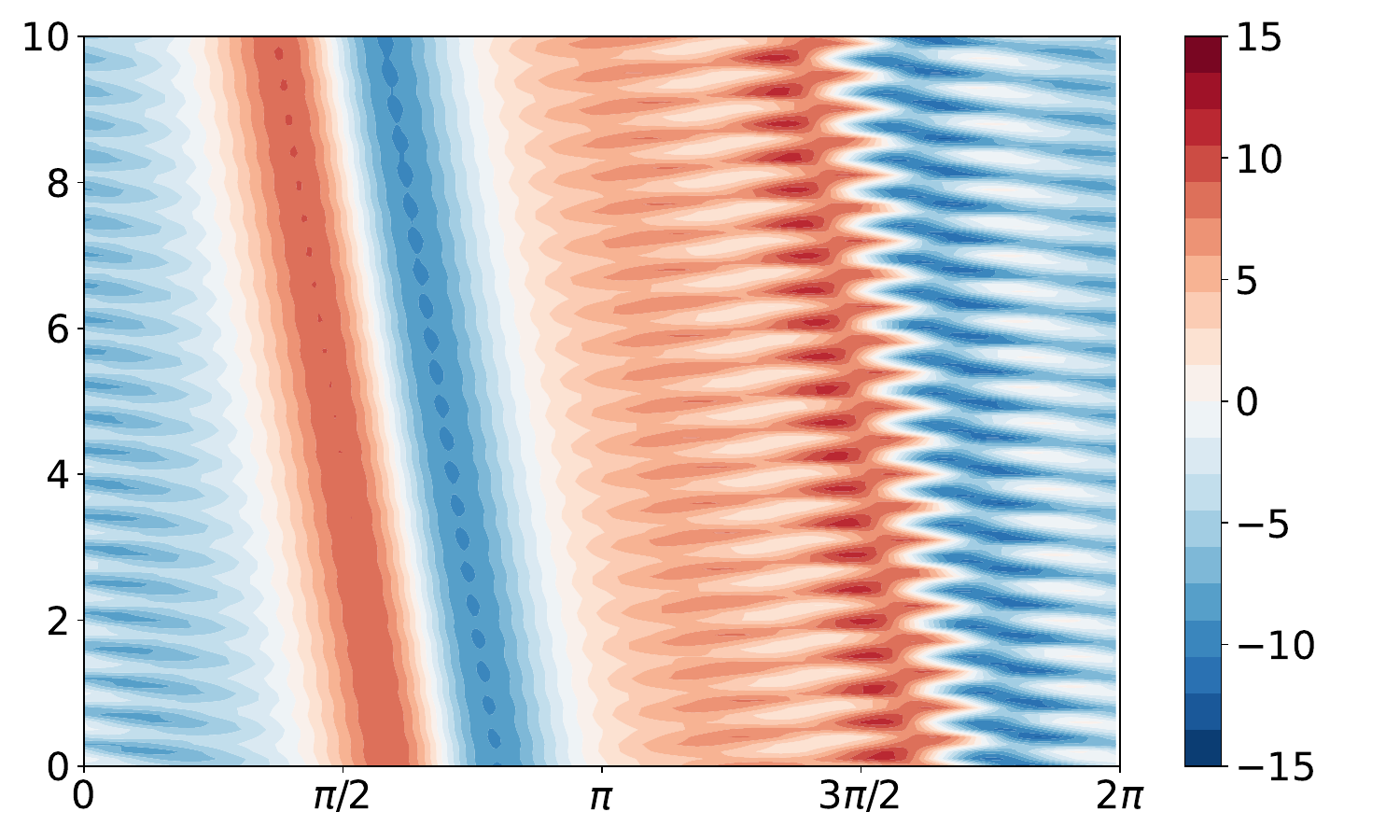}
      \node[anchor=south, align=center] at (0.45, 1.06) {\large S-R OpInf w/};
      \node[anchor=south, align=center] at (0.45, 0.95) {\large optimal penalties, 2.22\%};
      \node at (0.45, -0.05) {$x$};
      \node at (0.45, -0.18) {(d)};
    \end{tikzonimage}
    \captionsetup{labelformat=empty}
    \caption{\mbox{}}
    \vspace{-1.5em}
    \label{fig:ROM_training_snapshots_panel_d}
  \end{subfigure}
  \caption{Reconstructions of the base full-order solution by the four ROMs with $n = 5$ modes, each labeled with its relative error above the panel.
    The unregularized S-R OpInf ROM, panel~(b), diverges at $t\approx1.8$, resulting in an undefined relative error.
    See Fig.~\ref{fig:original_and_template_fitted_training_snapshots_panel_a}
    for comparison with the full-order solution snapshots.}
  \label{fig:ROM_training_snapshots}
\end{figure}

In this section, we discuss the reconstruction of the base traveling-beating solution.
Before turning to the ROMs, we examine the template fitting itself.
Fig.~\ref{fig:original_and_template_fitted_training_snapshots} shows the solution before (panel~a) and after (panel~b) template fitting.
Although the dominant drift is removed, a slight spatial fluctuation visible as a zig-zag pattern remains in the fitted profile.
This residual reflects a limitation of the chosen template-fitting approach, which tracks only global translations.
As a result, the secondary spatial movements of the beating wave are not entirely factored out by the continuous shift operator.

To quantify the reconstruction we use the relative error
\begin{equation}
    \label{eq:036_relative_error}
    \varepsilon = \sqrt{\frac{\sum_{m = 0}^{N_t-1} \|u_{\mathrm{ROM}}(t_m) - u(t_m)\|^2}{\sum_{m = 0}^{N_t-1} \|u(t_m)\|^2}},
\end{equation}
where $u_{\mathrm{ROM}}(t)$ is the full-order state reconstructed from the ROM trajectory.

Fig.~\ref{fig:ROM_training_snapshots} compares the reconstructions of the four
ROMs listed in Section~\ref{sec:numerics} with dimension $n = 5$.
Without regularization or re-projection, the S-R OpInf ROM (panel~b) diverges at $t\approx1.8$.
In contrast, the S-R OpInf ROM trained on the re-projected dataset (panel~c) is numerically stable during the online evaluation time window and reproduces the drift and modulation with a relative error of $6.12\%$, identical to that of the intrusive S-R Galerkin ROM (panel~a).
This confirms that re-projection drives the learned operators to their Galerkin counterparts (Appendix~\ref{app:reprojection}).
Nevertheless, activating penalization in the optimization problems of S-R OpInf turns out to be even more effective.
With the grid-searched optimal penalties, the S-R OpInf ROM (panel~d) attains a relative error of $2.22\%$, lower than that of the intrusive S-R Galerkin ROM.
These results demonstrate that with appropriate regularization, the non-intrusive S-R ROM can surpass the accuracy of the intrusive approach.
Nevertheless, it still depends on either regularization or re-projection to maintain numerical stability.

\begin{figure}[tbp]
  \centering
  \begin{tikzonimage}[width=0.62\linewidth]{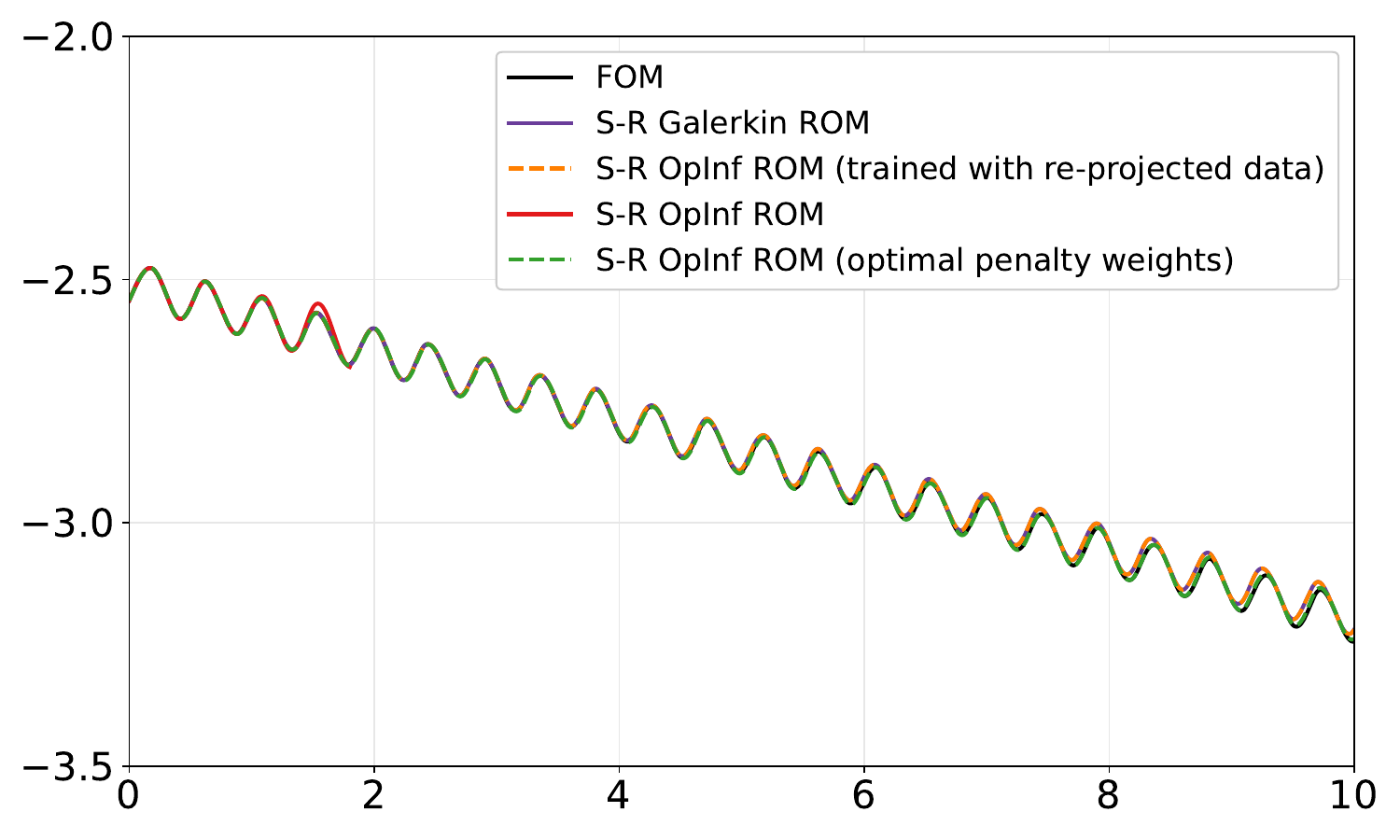}
    \node at (0.5, -0.06) {$t$};
    \node[overlay] at (-0.07, 0.5) {$c(t)$};
  \end{tikzonimage}
  \caption{The shift amount $c(t)$ of the base solution recovered by the FOM and
    the four $5$-dimensional ROMs.  The unregularized S-R OpInf ROM terminates at $t\approx1.8$; the remaining models track the FOM, with the optimal-penalty S-R OpInf ROM the most accurate.}
  \label{fig:shift_amount}
\end{figure}

Apart from comparing the full-order fields $u(x, t)$ generated by different models, we also evaluate the performance of these models in recovering the shift amount.
Fig.~\ref{fig:shift_amount} shows $c(t)$ recovered by each model.
All stable ROMs track the true shift amount accurately.
Among them, the S-R OpInf ROM with optimal penalties is the most accurate, while the S-R Galerkin and the re-projected S-R OpInf ROMs accumulate a slight drift accumulating over time.
Since a small misalignment in shift amount can produce a large discrepancy in $u(x, t)$ at the original frame, accurate recovery of $c(t)$ is essential.
Note that $\dot c(t)$ is neither constant nor sign-definite, as it reflects the oscillatory beating of the solution and the its continuous alignment with the template.

\subsection{Generalization across multiple trajectories}
\label{sec:generalization-multiple}

We now turn to the second task: reconstructing the family of perturbed solutions described in Section~\ref{sec:numerics}.
In this scenario, S-R ROMs based on a single POD basis and a single set of reduced operators are built to represent a diverse ensemble of trajectories, making this task substantially more challenging than the single-solution case.
Since the perturbations are linearly unstable, each trajectory undergoes an initial transient during which the perturbation is amplified, causing initially nearby trajectories to separate from one another over time.
Consequently, model performance is expected to degrade non-uniformly as the trajectories separate,
and a single time-averaged error metric is insufficient.
To explicitly track how the reconstruction quality evolves dynamically, we introduce the cumulative relative error at every snapshot time $t_k$:
\begin{equation}
    \label{eq:037_cumulative_relative_error}
    \varepsilon_{\text{c}}(t_k) = \sqrt{\frac{\sum_{m = 0}^{k} \|u_{\mathrm{ROM}}(t_m) - u(t_m)\|^2}{\sum_{m = 0}^{k} \|u(t_m)\|^2}},
\end{equation}
which is evaluated by truncating both sums in~\eqref{eq:036_relative_error}.
From its definition, we immediately have $\varepsilon_{\text{c}}(t_{N_t-1}) = \varepsilon$.

\begin{figure}[tbp]
  \centering
  \begin{subfigure}[b]{0.48\textwidth}
    \centering
    \begin{tikzonimage}[width=\linewidth]{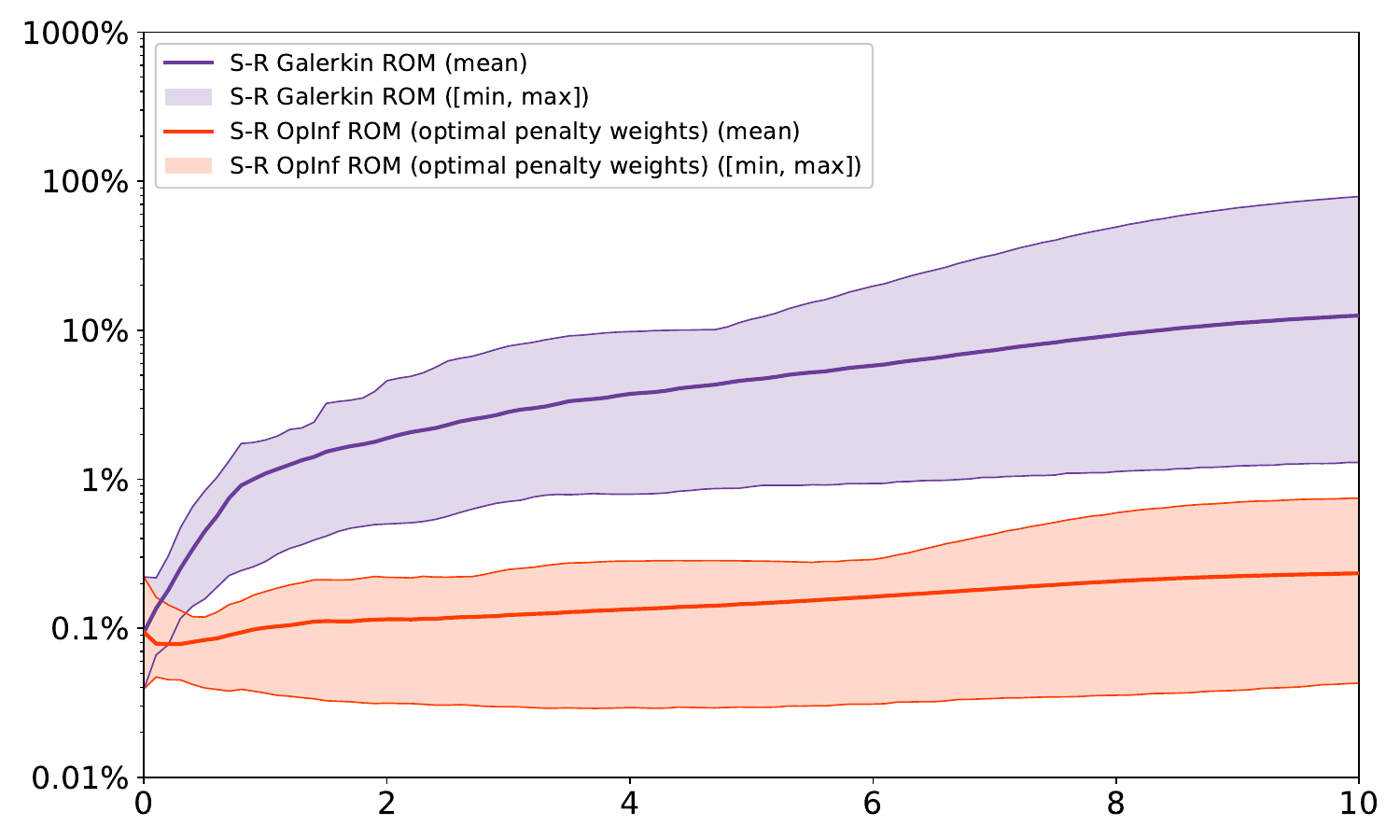}
      \node at (0.45, -0.05) {$t$};
      \node at (0.45, -0.18) {(a)};
    \end{tikzonimage}
    \captionsetup{labelformat=empty}
    \caption{\mbox{}}
    \vspace{-1.5em}
    \label{fig:cumulative_error_bands_panel_a}
  \end{subfigure}
  \hfill
  \begin{subfigure}[b]{0.48\textwidth}
    \centering
    \begin{tikzonimage}[width=\linewidth]{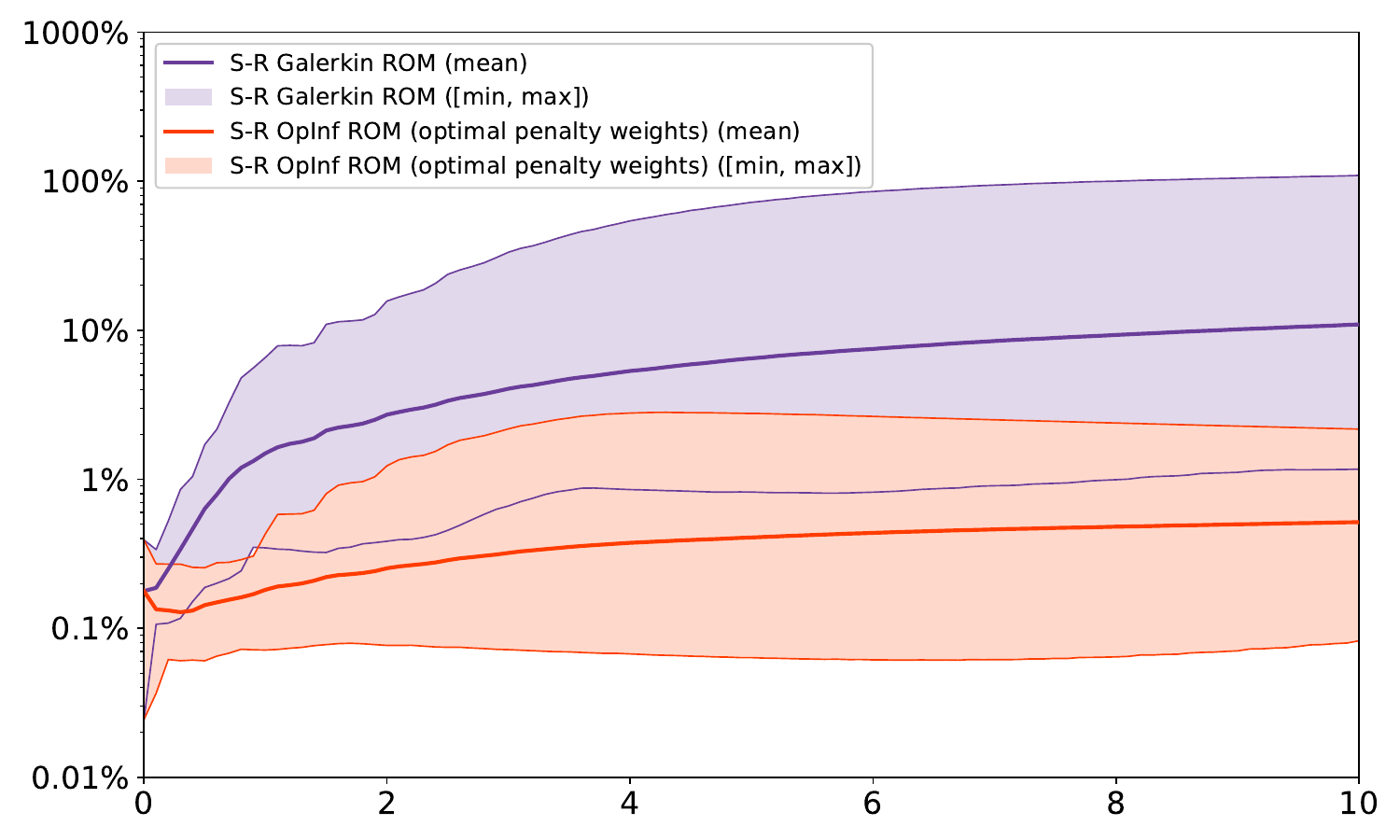}
      \node at (0.45, -0.05) {$t$};
      \node at (0.45, -0.18) {(b)};
    \end{tikzonimage}
    \captionsetup{labelformat=empty}
    \caption{\mbox{}}
    \vspace{-1.5em}
    \label{fig:cumulative_error_bands_panel_b}
  \end{subfigure}
  \caption{Cumulative relative error \eqref{eq:037_cumulative_relative_error} of the S-R Galerkin ROM (purple) and the S-R OpInf ROM with the optimal penalties (orange) with dimension $n = 20$.
    The trajectory-wise mean cumulative error is depicted in solid lines, while the minimal and maximal errors are the borders of the shaded regions.
    Panel (a) shows the result over the $10$ training trajectories at the perturbation intensity $\gamma = 10\%$, and panel (b) over the $30$ testing trajectories at $\gamma = 12\%$.
    No trajectory diverges for either model.}
  \label{fig:cumulative_error_bands}
\end{figure}

Figure~\ref{fig:cumulative_error_bands}(a) reports the mean cumulative relative
error and the trajectory-wise $[\min,\max]$ envelope over the $10$ training
trajectories.  Both the S-R Galerkin ROM and the penalty-optimized S-R OpInf ROM
built upon $n = 20$ POD modes remain numerically stable on all trajectories, but their error behavior is different.
The S-R OpInf error remains nearly flat over the integration window,
whereas the S-R Galerkin error grows more rapidly and develops a substantially
wider upper envelope.  At the final time $t=10$, the S-R OpInf ROM has a
trajectory-wise mean relative error of only $0.23\%$, with errors ranging from
$0.04\%$ to $0.70\%$.  In contrast, the S-R Galerkin ROM has a mean error of
$12.32\%$, with a much wider range $[1.31\%,76.90\%]$.
The discrepancy in the predictive performance of the two ROMs, observed both here and in the single-solution task, highlights the essential stabilizing role of the regularization procedure intrinsic to the Operator Inference framework~\cite{Swischuk2020, Jain2021, McQuarrie2021, Sawant2021, Kramer2024}.
Indeed, truncating a high-dimensional fluid system inherently discards high-frequency dissipative modes that are important for preventing unphysical energy growth, resulting in a potentially misspecified reduced model~\cite{Kramer2024}.
In contrast, S-R OpInf actively penalizes the spurious dynamics that arises from overfitted operators with excessively large norms.
When coupled with a grid-search procedure, this regularization serves as a crucial compensatory mechanism for the unresolved dynamics, which eventually yields a more robust model than straightforward Galerkin projection.

Although the S-R OpInf ROM performs better than the S-R Galerkin ROM on the training dataset, this advantage is expected since the S-R OpInf operators are learned from the same trajectories being reconstructed.
A more stringent test is whether the learned operators remain predictive for trajectories initialized from unseen perturbations.
Figure~\ref{fig:cumulative_error_bands}(b) repeats the comparison on the $30$ testing trajectories with perturbation intensity $\gamma=12\%$, which represents the nearest out-of-sample condition.
The performance of the S-R OpInf ROM is only mildly affected: its mean error increases slightly to $0.55\%$, with a trajectory-wise range of $[0.09\%, 2.36\%]$.
For comparison, the S-R Galerkin ROM yields a mean error of $11.43\%$, with its trajectory-wise maximal error extending to $116.68\%$.
While this testing mean is surprisingly lower than the corresponding training error ($12.32\%$),
this does not indicate superior generalization.
Instead, this fluctuation in trajectory-wise mean testing error, combined with an increased variance in testing errors across the ensemble, highlights the inherent sensitivity of the uncalibrated Galerkin operators to different initial realizations.
In contrast, the S-R OpInf ROM maintains consistent accuracy across both datasets. This confirms that the embedded regularization not only stabilizes the learned dynamics but mitigates overfitting to the training set.
Consequently, the data-driven model is endowed with robust generalization capabilities, thereby preserving its order-of-magnitude accuracy advantage even outside the training distribution.

\begin{figure}[tbp]
  \centering
  \begin{tikzonimage}[width=\linewidth]{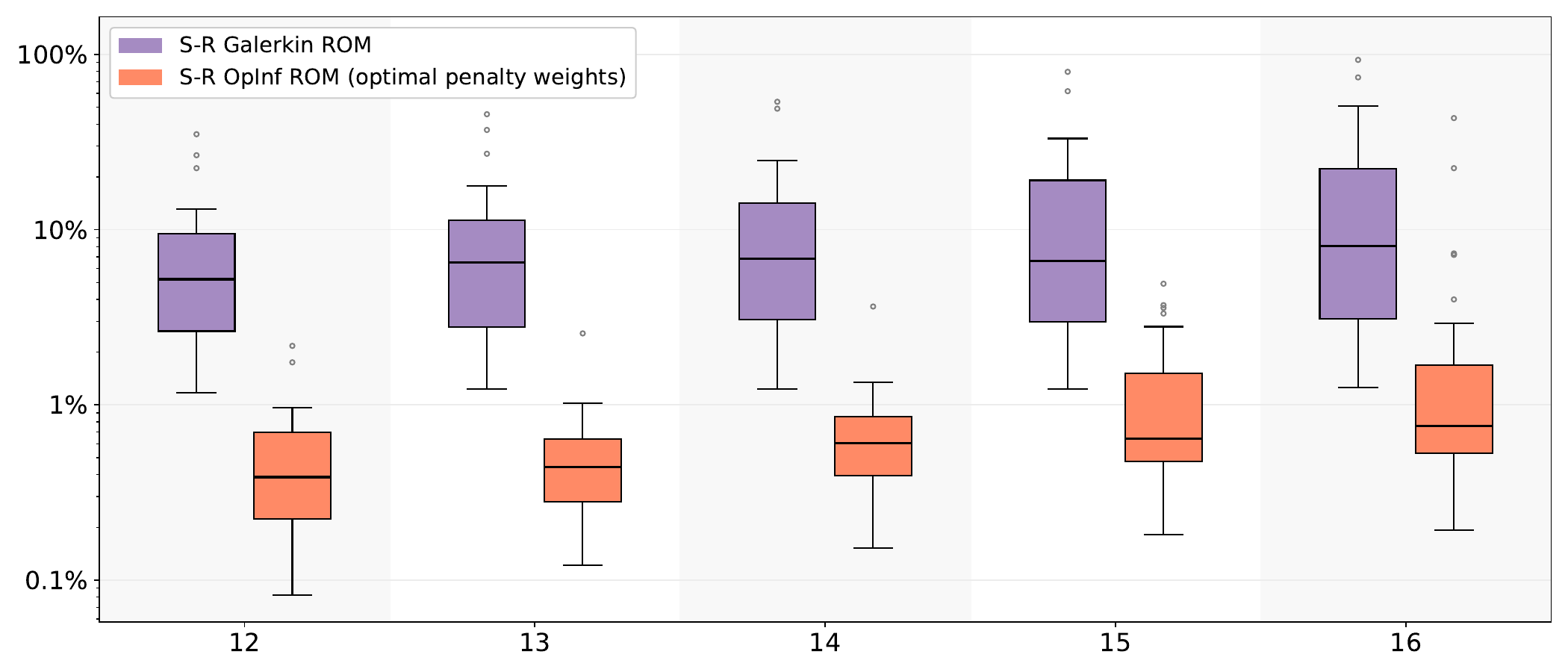}
    \node at (0.5, -0.07) {perturbation intensity $\gamma$ (\%)};
  \end{tikzonimage}
  \caption{Total relative error~\eqref{eq:036_relative_error} on the testing trajectories, grouped by perturbation intensity $\gamma$ ($30$ trajectories per group), for the S-R Galerkin ROM (purple) and the S-R OpInf ROM with the optimal penalties (orange).  Boxes span the interquartile range, with the median marked in black; whiskers extend to the most extreme values within $1.5$ interquartile ranges of the box, and circles mark values beyond.}
  \label{fig:perturbation_boxplot}
\end{figure}

To extend the analysis above, we evaluate both ROMs on testing trajectories generated at all selected perturbation intensities $\gamma \in \{12, 13, 14, 15, 16\}\%$.
For each intensity level, 30 distinct trajectories are reconstructed (the $\gamma=12\%$ case utilizes the results already obtained above).
The trajectory-wise distribution of errors depicted in Figure~\ref{fig:perturbation_boxplot}.
Across this entire range, the S-R OpInf ROM demonstrates better accuracy and robustness, characterized by a substantially lower mean and variance in relative error. Specifically, as $\gamma$ increases, its median error stays below $1\%$, while its interquartile range remains tightly bounded below $2\%$.
Except a few extreme outliers at $\gamma = 16\%$, the trajectory-wise errors for S-R OpInf consistently remain below the $5\%$ level.
The S-R Galerkin ROM, on the other hand, yields median errors from about $5\%$ at $\gamma = 12\%$ to $8\%$ at $\gamma = 16\%$ and degrades more rapidly.
In addition, at least $10\%$ (3 out of 30) of the reconstructed realizations in every tested ensemble suffer severe accuracy degradation with relative errors exceeding $20\%$, whereas the S-R OpInf ROM maintains a low, tightly clustered error distribution across the entire sweep.

In summary, these numerical experiments establish that the non-intrusive S-R OpInf framework systematically surpasses its intrusive S-R Galerkin counterpart across the full ensemble of perturbed solutions.
It yields better average accuracy, maintains a tighter error variance, and effectively mitigates the occurrence of severely degraded trajectory reconstructions.
Given that the perturbation intensity $\gamma$ effectively parameterizes a family of dynamic behaviors,
these findings further indicate the S-R OpInf ROM as a promising generalizable framework for parameter-dependent shift-equivariant systems, provided that a sufficient training dataset spans the parameter regime of interest.
We will revisit this prospect in Section~\ref{sec:conclusions}.

\section{Conclusion}
\label{sec:conclusions}

In this paper, we introduced a novel non-intrusive reduced-order modeling (ROM) framework for approximating the dynamics of shift-equivariant systems, which are commonly encountered in advection-dominated problems with spatially drifting solutions.
Our approach extends the standard Operator Inference (OpInf) methodology to construct data-driven ROMs, thereby entirely bypassing the need for individual full-order operators in the underlying governing equations.
To handle the shift equivariance, we incorporate symmetry-reduction techniques, including template fitting to isolate spatial translations in the solution, and the algebraic shift reconstruction equation to track the eliminated shift during online evaluation.
The resulting symmetry-reduced operator inference (S-R OpInf) ROM yields a low-order surrogate system for the efficient evaluation of both the fitted profile dynamics and the shift amount, enabling the full reconstruction of the original traveling field.

Training the S-R OpInf model involves minimizing the residual errors in both the polynomial dynamics and the numerator of the shift reconstruction equation.
As with standard OpInf algorithms, the resulting minimization problems are formulated as linear least-squares problems with Tikhonov regularization, which admit efficient closed-form solutions without the need for iterative optimization.
Our numerical results show that this regularization is vital for penalizing spurious growth and stabilizing the learned dynamics against unmodeled truncation errors.
Furthermore, we generalize the re-projection technique presented in~\cite{Peherstorfer2020} to our setting for continuous-time symmetry-reduced systems.
This generalization ensures that the learned data-driven S-R OpInf operators from the re-projected dataset are the exact operators derived via intrusive Galerkin projection in the S-R Galerkin ROM, given that the full-order right-hand side function is accessible and the data matrix in the least-squares problem has full column rank.

We evaluated the performance of this framework on the Kuramoto–Sivashinsky equation, a one-dimensional PDE exhibiting traveling beating-wave solutions.
In both single-trajectory and multiple-trajectory reconstruction tasks, the S-R OpInf ROM trained with grid-searched optimal penalty weights consistently delivers high-fidelity reconstructions of the full-order field.
The strength of our approach is highlighted in the multi-solution scenario: when evaluated against unseen testing trajectories, the S-R OpInf framework exhibited remarkable robustness.
It systematically outperformed its intrusive S-R Galerkin counterpart, actively suppressing the severe accuracy degradation that afflicted the Galerkin model.
Furthermore, we validated the effectiveness of the re-projection technique in reproducing projection-based reduced operators via data-driven learning.
When trained on re-projected data, the S-R OpInf ROM yields predictions that exactly match those of the intrusive S-R Galerkin model.
This consistency further confirms that our data-driven, non-intrusive framework is highly versatile: it not only achieves good accuracy and generalization ability but is also capable of emulating mathematically rigorous low-order operators from Galerkin projection.

A primary motivation for S-R OpInf is to model the low-rank dynamics of transport-dominated systems featuring shifting patterns.
For systems whose governing equations and boundary conditions are equivariant under translations along a periodic coordinate, this model can be used to efficiently track the time evolution of simple traveling waves and other relative equilibria.
Typical real-world examples are transitional and turbulent shear
flows, including pipe flow and channel flow where exact coherent structures such
as traveling waves and relative periodic orbits evolve through a combination of
shape dynamics and downstream or azimuthal drift~\cite{Wedin2004, Willis2013, Willis2016, Budanur2017}.
Similarly, in idealized annular reacting-flow models such as rotating detonation engines with azimuthally uniform injection, the dynamics are heavily dominated by wave motion along a periodic angular coordinate~\cite{Farcas2023}.
In a symmetry-reduced frame, these complex traveling structures become equilibria and relative periodic orbits become periodic or recurrent low-dimensional trajectories~\cite{Budanur2015, Marensi2023}.
By learning the frozen-frame dynamics directly from shifted snapshots and their velocities, S-R OpInf yields a compact, non-intrusive surrogate that overcomes the expensive cost of full-order simulations.
This computationally efficient framework directly enables downstream tasks that are otherwise cumbersome, including data assimilation, many-query design optimization, and the robust extraction of low-rank physical patterns.

While the proposed S-R OpInf framework demonstrates strong predictive performance
in the examples considered here, several methodological limitations delineate its
present scope. First, the current formulation uses a single global phase variable and
therefore a single co-moving frame. It can stationarize a dominant traveling
structure, or multiple structures sharing the same drift, but it cannot simultaneously
freeze components that propagate with different speeds or in different directions.
Multi-frame shifted decompositions and dynamically transformed mode formulations
provide possible starting points for this setting~\cite{Reiss2019, Black2020,
Burela2025}; extending these ideas to a fully non-intrusive operator-inference
framework remains an important direction for future work. Second, our use of the
term ``non-intrusive'' should be understood in the operator-learning sense.
The framework does not require access to individual full-order operators, but it still assumes access to high-fidelity snapshot data and to
full-state preprocessing operations needed to construct the frozen data and training matrices, including continuous shifts, derivative estimates, and inner products.
Last but not least, a complementary methodological direction is to combine symmetry reduction with non-intrusive constructions of oblique projection operators~\cite{Padovan2024}.
Such a combination could reduce the restriction imposed by orthogonal POD
projection and may improve the representation of transient dynamics that are not
well captured by the leading POD subspace alone.

Beyond the fixed-parameter, autonomous setting studied in this work, the S-R OpInf
framework also suggests two natural extensions.
The first is model reduction of parametric shift-equivariant systems.
Given training trajectories at several values of a PDE parameter, one may align the snapshots at each sampled parameter, assemble a common frozen POD basis, and infer a separate set of reduced operators for each parameter value.
Following parametric operator inference~\cite{Farcas2023, McQuarrie2023}, the operators associated with an unsampled parameter can then be obtained by interpolation, yielding a ROM that evolves in symmetry-reduced coordinates and reconstructs the corresponding physical trajectory through the shift dynamics.
A symmetry-specific issue is that the template used for fitting may itself vary with the parameter.
This can be handled either by using a parameter-dependent template or by selecting a single template that remains approximately co-aligned with the dominant coherent structure over the parameter range of interest.
The second extension is control-oriented reduced-order modeling via Operator Inference~\cite{Kramer2017}.
When the shift-equivariant full-order dynamics include an external input that is also shift-equivariant, the frozen-state equations and the shift reconstruction equation can be augmented with input-dependent terms.
The resulting S-R OpInf ROM would then provide a non-intrusive reduced model for
predicting responses to previously unseen input signals, subject to the usual
requirement that the training inputs sufficiently excite the relevant controlled
dynamics.

\backmatter

\bmhead{Supplementary information}

Not applicable.

\bmhead{Acknowledgements}

We thank the anonymous referees for many helpful suggestions, especially for
suggesting an improved setup for the minimization
problem~\eqref{eq:032_sr_opinf_loss}.

\section*{Declarations}

\subsection*{Funding}

We gratefully acknowledge support from the Air Force Office of Scientific Research, award FA9550-19-1-0005.

\subsection*{Conflict of Interest}

The authors declared that they have no conflict of interest.

\subsection*{Ethics approval and consent to participate}

Not applicable.

\subsection*{Consent for publication}

Not applicable.

\subsection*{Data availability}

All datasets and figures supporting the findings of this article are available from the corresponding author upon reasonable request.

\subsection*{Materials availability}

Not applicable.

\subsection*{Code availability}

The code used in this study is available at our public GitHub repository, \url{https://github.com/Yu-Shuai-PU/sr-opinf}, and is archived at Zenodo with DOI \url{https://doi.org/10.5281/zenodo.20823370}.  The numerical experiments were run with Python~3.11 on macOS~14.x; the pinned dependency versions are listed in \texttt{requirements.txt} in the repository.

\subsection*{Author Contributions}

Yu Shuai and Clarence W. Rowley contributed to the study conception.
Algorithm design, material preparation, data collection and analysis were performed by Yu Shuai.
The first draft of the manuscript was written by Yu Shuai.
Yu Shuai and Clarence W. Rowley worked on revisions to the manuscript.
Yu Shuai and Clarence W. Rowley read and approved the final manuscript.

\begin{appendices}
\section{Re-projection and exact recovery of the symmetry-reduced Galerkin ROM}
\label{app:reprojection}
\renewcommand{\theequation}{A.\arabic{equation}}
\renewcommand{\theHequation}{A.\arabic{equation}}
\setcounter{equation}{0}
\renewcommand{\thetheorem}{A.\arabic{theorem}}
\renewcommand{\theHtheorem}{A.\arabic{theorem}}
\setcounter{theorem}{0}
\renewcommand{\theremark}{A.\arabic{remark}}
\renewcommand{\theHremark}{A.\arabic{remark}}
\setcounter{remark}{0}

This appendix establishes the result used in Section~\ref{sec:numerics}: when symmetry-reduced operator inference is trained on a re-projected dataset, its learned operators coincide \emph{exactly} with those of the symmetry-reduced Galerkin ROM.  We build the result in three steps.  We first recall the discrete-time re-projection technique of~\cite{Peherstorfer2020} in our function-space notation (Appendix~\ref{app:discrete}); we then extend it to continuous-time dynamical systems, which to our knowledge is new and rests on a single assumption---direct access to the full-order velocity field as a whole (Appendix~\ref{app:continuous}); and we finally combine that extension with symmetry reduction (Appendix~\ref{app:sr}).
In solving the OpInf optimization problem, we use the duplicate-free convention introduced in \ref{sec:opinf}.

Throughout the appendix, $u_0\in\U$ denotes the template function of the main text, $\partial_x$ the spatial-derivative operator (skew-adjoint on the periodic domain, $\ip{\partial_x v}{w}=-\ip{v}{\partial_x w}$), and $P v = \sum_{i=1}^{n}\ip{v}{\varphi_i}\varphi_i$ the orthogonal projection onto the reduced subspace $ \U_n=\mathrm{span}\{\varphi_i\}_{i=1}^n$.
Learned coefficients carry a hat ($\hat A,\hat B,\hat p,\hat Q$) and Galerkin-projected ones a tilde ($\tilde A,\tilde B,\tilde p,\tilde Q$).

\subsection{Discrete-time re-projection}
\label{app:discrete}

Consider an autonomous discrete-time system with quadratic nonlinearity,
\begin{equation}
\label{eq:app-fom-discrete}
u_{m+1} = f(u_m) = A\,u_m + B(u_m,u_m),
\qquad m = 0,\dots,K-1,
\end{equation}
with full state $u_m\in\U$, a linear operator $A:\U\to\U$, and a symmetric bilinear operator $B:\U\times\U\to\U$, $B(v,w)=B(w,v)$.  Write the reduced coordinates as $a_{m,i} = \ip{u_m}{\varphi_i}$, collected into $a_m\in\mathbb{R}^n$.  There are two routes to a reduced model.

If the operators are available, intrusive Galerkin projection forms
\begin{equation}
\label{eq:app-galerkin-discrete}
\tilde{A}_{ij} = \ip{A\varphi_j}{\varphi_i},
\qquad
\tilde{B}_{ijk} = \ip{B(\varphi_j,\varphi_k)}{\varphi_i}
\qquad (\tilde B_{ijk}=\tilde B_{ikj}),
\end{equation}
and yields the Galerkin reduced model
\begin{equation}
\label{eq:app-rom-discrete}
\tilde{a}_{m+1,i} = \ip{f\!\Big(\textstyle\sum_{j}\tilde a_{m,j}\varphi_j\Big)}{\varphi_i}
= \tilde{A}_{ij}\,\tilde{a}_{m,j} + \tilde{B}_{ijk}\,\tilde{a}_{m,j}\tilde{a}_{m,k},
\qquad \tilde{a}_{0,i} = \ip{u_{\mathrm{init}}}{\varphi_i},
\end{equation}
where $u_{\mathrm{init}}$ is the initial full state.  If the operators are unavailable, operator inference instead fits $\hat{A},\hat{B}$ to the projected snapshots by least squares,
\begin{equation}
\label{eq:app-opinf-discrete}
\min_{\hat{A},\hat{B}}\;
\sum_{m=0}^{K-1}\sum_{i=1}^{n}\Big[\, \hat{A}_{ij}\,a_{m,j} + \hat{B}_{ijk}\,a_{m,j}a_{m,k} - a_{m+1,i}\,\Big]^2 .
\end{equation}

In general, solving~\eqref{eq:app-opinf-discrete} with projected FOM
snapshot pairs does not recover the Galerkin reduced operators. To see this, write
\[
    u_m = Pu_m + r_m = \sum_{i}a_{m, i}\varphi_i + r_m, \qquad Pu_m \in \U_n,\quad r_m \perp \U_n,
  \]
from which the next projected snapshot is generated as:
\[
    a_{m+1, i}
    = \ip{f(u_m)}{\varphi_i}
    = \ip{f\bigg(\sum_{i}a_{m, i}\varphi_i + r_m\bigg)}{\varphi_i}.
\]
Thus $a_{m+1}$ generally depends on the unresolved component $r_m$ and is not a
function of $a_m$ alone.
Equivalently, the exact projected dynamics are not closed,
or Markovian, in the reduced coordinates. By contrast, a Galerkin reduced
model~\eqref{eq:app-rom-discrete} imposes a Markovian closure by evolving only
states in $\U_n$. Consequently, the operators $\hat A,\hat B$ inferred from projected FOM snapshot
pairs act as best-fit Markovian closure operators for the observed projected
trajectory.
They generally differ from the Galerkin operators $\tilde A,\tilde B$.
This distinction is closely related to the memory and closure terms that arise in the
Mori--Zwanzig formalism, see~\cite{Peherstorfer2020} for further reference.

The remedy is to generate training data that are Markovian by construction. Starting from the initial condition $\bar{a}_{0,i} = \ip{u_{\mathrm{init}}}{\varphi_i} = \tilde{a}_{0, i}$, at each step $m$, we first lift the reduced state $\bar{a}_m$ back to the full-dimensional space via $\bar{a}_m \mapsto \sum_j \bar{a}_{m,j}\varphi_j$, which we refer to as the \textit{re-projected} full state.
Rather than advancing the true full-order state $u_m$, we advance this re-projected state by one full-order step and immediately project it back onto the basis:
\begin{equation}
\label{eq:app-reproj-discrete}
\bar{a}_{m+1,i} = \ip{f\!\Big(\textstyle\sum_{j}\bar a_{m,j}\varphi_j\Big)}{\varphi_i},
\qquad m = 0,\dots,K-1.
\end{equation}
This ensures that any off-subspace component is strictly discarded before it can feed into the subsequent step.
The sequence $\{\bar{a}_m\}_{m = 0}^{K-1}$ generated by this iterative procedure forms the \textit{re-projected} reduced coordinates. We next show that applying the OpInf learning procedure to this dataset analytically reproduces the projection-based Galerkin reduced operators.

\begin{proposition}[Re-projected trajectory equals the Galerkin trajectory]
\label{prop:app-reproj-discrete}
Let $\{\bar{a}_m\}_{m=0}^{K-1}$ be generated by~\eqref{eq:app-reproj-discrete} and $\{\tilde{a}_m\}_{m=0}^{K-1}$ by the Galerkin model~\eqref{eq:app-rom-discrete}, both from the same initial condition $\bar{a}_{0}=\tilde{a}_{0}$.  Then $\bar{a}_m=\tilde{a}_m$ for all $m$; in particular the re-projected trajectory has zero closure error.
\end{proposition}

\begin{proof}
The key algebraic fact is the bilinearity of $B$: for any $a\in\mathbb{R}^n$,
\begin{equation}
\label{eq:app-bilinear}
B\!\Big(\textstyle\sum_{j}a_j\varphi_j,\ \sum_{k}a_k\varphi_k\Big)
= \sum_{j,k} a_j a_k\, B(\varphi_j,\varphi_k).
\end{equation}
Combined with the definitions~\eqref{eq:app-galerkin-discrete}, the re-projection map evaluated at any reduced state $a\in\mathbb{R}^n$ satisfies, for each $i$,
\begin{equation}
\label{eq:app-identity-discrete}
\ip{f\!\Big(\textstyle\sum_{j}a_j\varphi_j\Big)}{\varphi_i}
= \sum_{j} a_j\,\ip{A\varphi_j}{\varphi_i}
 + \sum_{j,k} a_j a_k\,\ip{B(\varphi_j,\varphi_k)}{\varphi_i}
= \tilde{A}_{ij}\,a_j + \tilde{B}_{ijk}\,a_j a_k ,
\end{equation}
using the linearity of $A$ and of $\ip{\cdot}{\varphi_i}$, and the bilinearity~\eqref{eq:app-bilinear}.  That is, the re-projection map coincides identically with the Galerkin right-hand side~\eqref{eq:app-rom-discrete}.  We argue by induction.  For $m=0$ the claim holds by hypothesis.  Assuming $\bar{a}_m=\tilde{a}_m$, we have
\[
\bar{a}_{m+1,i}
\overset{\eqref{eq:app-reproj-discrete}}{=} \ip{f\!\Big(\textstyle\sum_{j}\bar a_{m,j}\varphi_j\Big)}{\varphi_i}
\overset{\eqref{eq:app-identity-discrete}}{=} \tilde{A}_{ij}\,\bar{a}_{m,j} + \tilde{B}_{ijk}\,\bar{a}_{m,j}\bar{a}_{m,k}
= \tilde{a}_{m+1,i}.
\]
Hence $\bar{a}_m=\tilde{a}_m$ for all $m$.
\end{proof}

Proposition~\ref{prop:app-reproj-discrete} says the re-projected data lie exactly on a Galerkin trajectory, so feeding them to operator inference recovers the Galerkin coefficients.

\begin{corollary}[Exact recovery, discrete time]
\label{cor:app-discrete}
Suppose the re-projected pairs $\{(\bar{a}_m,\bar{a}_{m+1})\}$ are rich enough that the regression matrix of~\eqref{eq:app-opinf-discrete}, whose rows are the feature vectors $[\,a_{m,j},\ a_{m,j}a_{m,k}\,]$, has full column rank, which requires at least
\begin{equation}
\label{eq:app-rank}
K \;\ge\; n + \frac{n(n+1)}{2}
\end{equation}
re-projected samples---the $n$ linear coefficients plus the $n(n+1)/2$ distinct quadratic monomials $a_ja_k$ (duplicates $a_ja_k=a_ka_j$ identified).  Then operator inference~\eqref{eq:app-opinf-discrete}, applied to the re-projected data, has a unique minimizer attaining zero residual, equal to the Galerkin coefficients: $\hat{A}=\tilde{A}$ and $\hat{B}=\tilde{B}$.
\end{corollary}

The price of re-projection is that the full-order map $f$ must be queryable at the selected states $\sum_i\bar a_{m,i}\varphi_i$, rather than only along a single given trajectory.

\subsection{Extension to continuous-time dynamical systems}
\label{app:continuous}

We now generalize the formulation in~\cite{Peherstorfer2020} to dynamical systems continuous in time.
The full-order model is now the autonomous quadratic evolution equation in $\U$,
\begin{equation}
\label{eq:app-fom-cont}
\partial_t u(t) = f(u(t)) = A\,u(t) + B\big(u(t),u(t)\big),
\end{equation}
with the same $A,B$ and basis $\{\varphi_i\}$, and reduced coordinates $a_i(t) = \ip{u(t)}{\varphi_i}$.  Galerkin projection uses the same coefficients~\eqref{eq:app-galerkin-discrete} and gives the reduced model $\dot{\tilde a}_i = \tilde A_{ij}\tilde a_j + \tilde B_{ijk}\tilde a_j\tilde a_k$.  Continuous-time operator inference fits $\hat A,\hat B$ to pairs of projected states and projected \emph{velocities}, $a_i(t_m) = \ip{u(t_m)}{\varphi_i}$ and $\dot a_i(t_m) = \ip{f(u(t_m))}{\varphi_i}$, by
\begin{equation}
\label{eq:app-opinf-cont}
\min_{\hat{A},\hat{B}}\;
\sum_{m=0}^{K-1}\sum_{i=1}^{n}\Big[\, \hat{A}_{ij}\,a_j(t_m) + \hat{B}_{ijk}\,a_j(t_m)a_k(t_m) - \dot{a}_i(t_m)\,\Big]^2 .
\end{equation}
As in the case of discrete-time systems, solving this minimization problem does not return the Galerkin coefficients.
The snapshot can be similarly decomposed as $u(t_m) = P u(t_m) + (I-P)u(t_m)$, so the target velocity $\dot a_i(t_m)=\ip{f(u(t_m))}{\varphi_i}$ is polluted by the off-subspace component $(I-P)u(t_m)$, which generally makes the learned $\hat A,\hat B$ differ from $\tilde A,\tilde B$.

The remedy is to remove the off-subspace component \emph{before} querying the vector field: for each snapshot we project onto $\U_n$ and evaluate the full-order right-hand side at the lifted state:
\begin{equation}
\label{eq:app-reproj-cont}
\bar{a}_i(t_m) = \ip{u(t_m)}{\varphi_i},
\qquad
\bar{v}_i(t_m) = \ip{f\big(P u(t_m)\big)}{\varphi_i}
= \ip{f\!\Big(\textstyle\sum_{j}\bar a_j(t_m)\varphi_j\Big)}{\varphi_i},
\end{equation}
and use the re-projected pairs $\{(\bar{a}(t_m),\bar{v}(t_m))\}$ in place of $\{(a(t_m),\dot{a}(t_m))\}$.  We make explicit the conditions (or assumptions) under which this recovers the Galerkin model:
\begin{itemize}
\item[\textbf{(A1)}] \textbf{Direct vector-field access.}  The full-order right-hand side function $f$ can be queried as a whole at selected states: for each $\bar{a}(t_m)$, the value $f(\sum_j\bar a_j(t_m)\varphi_j)$---and hence $\bar{v}_i(t_m)$---is evaluated \emph{exactly}, not approximated by a finite difference of snapshots.
\item[\textbf{(A2)}] \textbf{Quadratic full-order model.}  $f$ has the form~\eqref{eq:app-fom-cont}.
\item[\textbf{(A3)}] \textbf{Orthonormal reduced basis.}  $\ip{\varphi_i}{\varphi_j}=\delta_{ij}$.
\item[\textbf{(A4)}] \textbf{Data richness.}  The re-projected states $\{\bar{a}(t_m)\}$ make the regression matrix of~\eqref{eq:app-opinf-cont} of full column rank; as in the discrete case this requires $K \ge n + n(n+1)/2$ samples, with duplicate quadratic monomials identified.
\end{itemize}

\begin{proposition}[Exact matching of the re-projected vector field]
\label{prop:app-cont}
Under \textnormal{(A2)--(A3)}, for every reduced state $\bar{a}\in\mathbb{R}^n$,
\begin{equation}
\label{eq:app-identity-cont}
\ip{f\!\Big(\textstyle\sum_{j}\bar a_j\varphi_j\Big)}{\varphi_i}
= \tilde{A}_{ij}\,\bar{a}_j + \tilde{B}_{ijk}\,\bar{a}_j\bar{a}_k.
\end{equation}
\end{proposition}

\begin{proof}
Identity~\eqref{eq:app-identity-cont} is the algebraic identity~\eqref{eq:app-identity-discrete}: by the bilinearity~\eqref{eq:app-bilinear} of $B$, the linearity of $A$, and the definitions~\eqref{eq:app-galerkin-discrete},
\[
\ip{f\!\Big(\textstyle\sum_{j}\bar a_j\varphi_j\Big)}{\varphi_i}
= \sum_{j}\bar a_j\,\ip{A\varphi_j}{\varphi_i}
 + \sum_{j,k}\bar a_j\bar a_k\,\ip{B(\varphi_j,\varphi_k)}{\varphi_i}
= \tilde{A}_{ij}\,\bar{a}_j + \tilde{B}_{ijk}\,\bar{a}_j\bar{a}_k . \qedhere
\]
\end{proof}

\begin{remark}[No time integration or full-order trajectory is required]
\label{rem:app-no-induction}
The discrete-time Proposition~\ref{prop:app-reproj-discrete} needed an induction over time steps, because there the re-projected states are linked by the one-step map.  In continuous time the training data are instantaneous (state, velocity) pairs and the matching identity~\eqref{eq:app-identity-cont} holds pointwise at each sampled state independently; no induction is needed.  The re-projected states $\{\bar{a}(t_m)\}$ therefore need not form a trajectory at all---any collection rich enough to satisfy \textnormal{(A4)} suffices.  In particular one may \emph{synthesize} the dataset by choosing any $K\ge n+n(n+1)/2$ reduced states in general position, lifting each to $\sum_j\bar a_j\varphi_j$, and evaluating $\ip{f(\cdot)}{\varphi_i}$, decoupling the training data entirely from trajectories of the full-order model.
\end{remark}

\begin{corollary}[Exact recovery in continuous time]
\label{cor:app-cont}
Under \textnormal{(A1)--(A4)}, the operator inference problem~\eqref{eq:app-opinf-cont} applied to the re-projected data $\{(\bar{a}(t_m),\bar{v}(t_m))\}$ of~\eqref{eq:app-reproj-cont} has a unique minimizer attaining zero residual, equal to the Galerkin coefficients: $\hat{A}=\tilde{A}$ and $\hat{B}=\tilde{B}$.
\end{corollary}

\begin{proof}
By Proposition~\ref{prop:app-cont} every re-projected pair satisfies $\bar{v}_i(t_m) = \tilde{A}_{ij}\,\bar{a}_j(t_m) + \tilde{B}_{ijk}\,\bar{a}_j(t_m)\bar{a}_k(t_m)$, so $(\tilde{A},\tilde{B})$ drives the nonnegative loss~\eqref{eq:app-opinf-cont} to zero and is a global minimizer.  The problem is a linear least-squares fit; when the design matrix of \textnormal{(A4)} has full column rank the minimizer is unique, the duplicate quadratic monomials identified as in Corollary~\ref{cor:app-discrete}.  Hence $(\hat{A},\hat{B})=(\tilde{A},\tilde{B})$.
\end{proof}

\begin{remark}[Why a finite difference is not enough]
\label{rem:app-no-fd}
Assumption (A1) is what makes the recovery exact.  A finite difference of the \emph{recorded} trajectory, $(u(t_m+\Delta t)-u(t_m))/\Delta t\approx f(u(t_m))$, approximates the velocity at the actual snapshot $u(t_m)$, which still carries the off-subspace component $(I-P)u(t_m)$; it returns the polluted target, not the re-projected velocity.  The assumption can be weakened to a time-stepper re-initialized at selected states: integrating~\eqref{eq:app-fom-cont} for a short time $\Delta t$ from the lifted state $P u(t_m)$ and finite-differencing yields $\ip{f(Pu(t_m))}{\varphi_i}+O(\Delta t)$, recovering the Galerkin coefficients only up to $O(\Delta t)$.  Direct evaluation of $f$ at the lifted state removes this error entirely.
\end{remark}

\subsection{Symmetry-reduced re-projection}
\label{app:sr}

We finally combine the continuous-time result with symmetry reduction.  Following the method of slices, the solution is written as a spatially shifted profile $u = S_{c}\,u_{\mathrm{fit}}$, with $S_c$ the shift operator, $c$ the shift amount, and $u_{\mathrm{fit}}\in\U$ the template-fitted profile.
Denote $\{\varphi_i\}$ as the POD basis computed from the template-fitted snapshots.
The reduced coordinates are computed as $a_i = \ip{u_{\mathrm{fit}}}{\varphi_i}$, and the intrusive symmetry-reduced Galerkin ROM for $a(t)$ and $c(t)$ are given as
\begin{equation}
\label{eq:app-sr-rom}
\dot{a}_i
= \tilde{A}_{ij}\,a_j + \tilde{B}_{ijk}\,a_j a_k + \dot c\,\tilde{C}_{ij}\,a_j,
\qquad
\dot c = -\frac{\tilde{p}_j a_j + \tilde{Q}_{jk}\,a_j a_k}{s_j a_j},
\end{equation}
with $\tilde A,\tilde B$ as in~\eqref{eq:app-galerkin-discrete}, the shift-coupling coefficients $\tilde{C}_{ij} = \ip{\partial_x\varphi_j}{\varphi_i}$, the denominator coefficients $s_j = \ip{\partial_x \varphi_j}{\partial_{x}u_0}$, and the numerator coefficients defined through
\begin{equation}
\label{eq:app-sr-pq}
\tilde{p}_j a_j + \tilde{Q}_{jk}\,a_j a_k
= \ip{f\!\Big(\textstyle\sum_{j}a_j\varphi_j\Big)}{\partial_x u_0},
\qquad
\tilde{p}_j = \ip{A\varphi_j}{\partial_x u_0},
\quad
\tilde{Q}_{jk} = \ip{B(\varphi_j,\varphi_k)}{\partial_x u_0} .
\end{equation}
The coefficients $\tilde C_{ij}$ and $s_j$ depend only on the basis and the template, so they need not be learned; only the $f$-dependent coefficients $\tilde A,\tilde B,\tilde p,\tilde Q$ carry information about the dynamics.  Without access to $A,B$, symmetry-reduced operator inference learns them by two least-squares problems,
\begin{align}
\label{eq:app-sr-ls-f}
\min_{\hat{A},\hat{B}} &\;
\sum_{m=0}^{K-1}\sum_{i=1}^{n}\Big[\, \hat{A}_{ij}\,a_j(t_m) + \hat{B}_{ijk}\,a_j(t_m)a_k(t_m)
- \ip{f(u_{\mathrm{fit}}(t_m))}{\varphi_i}\,\Big]^2, \\
\label{eq:app-sr-ls-c}
\min_{\hat{p},\hat{Q}} &\;
\sum_{m=0}^{K-1}\Big[\, \hat{p}_j\,a_j(t_m) + \hat{Q}_{jk}\,a_j(t_m)a_k(t_m)
- \ip{f(u_{\mathrm{fit}}(t_m))}{\partial_x u_0}\,\Big]^2 ,
\end{align}
which are~\eqref{eq:032a_poly_loss}--\eqref{eq:032b_num_loss} of the main text without regularization.
Both targets are evaluated at the actual fitted snapshots $u_{\mathrm{fit}}(t_m)$, which carry an off-subspace component, so as in Appendix~\ref{app:continuous} both are polluted and the learned coefficients differ from the symmetry-reduced Galerkin ones.

Re-projection cures both at once.  We project each fitted snapshot and query the full-order vector field at the lifted state, forming the two re-projected targets
\begin{equation}
\label{eq:app-sr-reproj}
\bar{a}_i(t_m) = \ip{u_{\mathrm{fit}}(t_m)}{\varphi_i},
\qquad
\bar{v}_i(t_m) = \ip{f\big(P u_{\mathrm{fit}}(t_m)\big)}{\varphi_i},
\qquad
\bar{g}(t_m) = \ip{f\big(P u_{\mathrm{fit}}(t_m)\big)}{\partial_x u_0},
\end{equation}
and use $\bar{v}(t_m),\bar{g}(t_m)$ in place of the polluted targets $\langle f(u_{\mathrm{fit}}(t_m)), \varphi_i\rangle$ and $\langle f(u_{\mathrm{fit}}(t_m)), \partial_x u_0\rangle$. Both require evaluating $f$ at the \emph{same} lifted state $P u_{\mathrm{fit}}(t_m)$---assumption \textnormal{(A1)} of Appendix~\ref{app:continuous}; no new assumption is introduced. We summarize symmetry-reduced operator inference with the re-projection technique in Algorithm~\ref{alg:sr-opinf-re-projection}.

\begin{algorithm}
  \caption{Symmetry-Reduced Operator Inference with Re-Projection}\label{alg:sr-opinf-re-projection}
  \begin{algorithmic}[1]
    \Require Snapshots of the solution ${u}(t_m)\in\U$. The template function $u_0\in\mathcal{U}$ for template fitting. The queryable full-order vector field $f(\cdot)$. The dimension of the desired
    reduced-order model $n$.
    \Ensure The Symmetry-Reduced Operator Inference model~\eqref{eq:033_sr_opinf_rom}, \eqref{eq:028_sr_opinf_shift_model} with coefficients identical to the projection-based (symmetry-reduced Galerkin) ones~\eqref{eq:010_galerkin_coeffs}, ~\eqref{eq:022_shift_coeffs}.
    \For{$m = 1, \dots, N_t$}
      \State Perform template fitting to obtain the shifted snapshots $u_{\mathrm{fit}}(t_m)$ according to~\eqref{eq:013_shifted_representation} and~\eqref{eq:015_template_fitting}.
    \EndFor
    \State Perform POD on the set of snapshots $\{u_{\mathrm{fit}}(t_m)\}$,
    to obtain orthonormal basis functions $\varphi_i\in\U$, $i=1,\ldots,n$.
    \For{$m = 1, \ldots, N_t$}
    \State Determine $\bar{a}_i(t_m)$, $\bar{v}_i(t_m)$ and $\bar{g}(t_m)$ from~\eqref{eq:app-sr-reproj}.
    \State Replace $a_i(t_m)$, $\langle f(u_{\mathrm{fit}}(t_m)), \varphi_i\rangle$ and $\langle f(u_{\mathrm{fit}}(t_m)), \partial_x u_0\rangle$ in~\eqref{eq:app-sr-ls-f}, \eqref{eq:app-sr-ls-c} with $\bar{a}_i(t_m)$, $\bar{v}_i(t_m)$ and $\bar{g}(t_m)$, respectively.
    \EndFor
    \State Determine the coefficients $s_j$ from~\eqref{eq:022_shift_coeffs}
    and $C_{ij}$ from~\eqref{eq:020_recon_coeffs}.
    \State Learn the coefficients $\hat{A}_{ij}$, $\hat{B}_{ijk}$, $\hat{p}_j$, and
    $\hat{Q}_{jk}$ by solving the optimization problems~\eqref{eq:app-sr-ls-f} and~\eqref{eq:app-sr-ls-c} with re-projected data $\bar{a}_i(t_m)$, $\bar{v}_i(t_m)$ and $\bar{g}(t_m)$ (noting the
    symmetry, $B_{ijk}=B_{ikj}$ and $Q_{jk}=Q_{kj}$).
    \State The reduced-order model is given by~\eqref{eq:033_sr_opinf_rom} with $\dot c$ given by~\eqref{eq:028_sr_opinf_shift_model}.
  \end{algorithmic}
\end{algorithm}

\begin{proposition}[Exact matching of both re-projected regressions]
\label{prop:app-sr}
Under \textnormal{(A2)--(A3)}, for every reduced state $\bar{a}\in\mathbb{R}^n$,
\begin{equation}
\label{eq:app-sr-identity}
\ip{f\bigg(\sum_j\bar a_j\varphi_j\bigg)}{\varphi_i} = \tilde{A}_{ij}\,\bar{a}_j + \tilde{B}_{ijk}\,\bar{a}_j\bar{a}_k,
\qquad
\ip{f\bigg(\sum_j\bar a_j\varphi_j\bigg)}{\partial_x u_0}
= \tilde{p}_j\,\bar{a}_j + \tilde{Q}_{jk}\,\bar{a}_j\bar{a}_k,
\end{equation}
with $\tilde{A},\tilde{B},\tilde{p},\tilde{Q}$ the symmetry-reduced Galerkin coefficients~\eqref{eq:app-galerkin-discrete},~\eqref{eq:app-sr-pq}.  Both re-projected targets lie exactly on the corresponding symmetry-reduced Galerkin right-hand sides; both regressions have zero closure error.
\end{proposition}

\begin{proof}
The first identity is Proposition~\ref{prop:app-cont}.  For the second, use the linearity of $\ip{\cdot}{\partial_x u_0}$ with the bilinearity~\eqref{eq:app-bilinear} of $B$:
\[
\ip{f\bigg(\sum_j\bar a_j\varphi_j\bigg)}{\partial_x u_0}
= \sum_{j}\bar a_j\,\ip{A\varphi_j}{\partial_x u_0}
+ \sum_{j,k}\bar a_j\bar a_k\,\ip{B(\varphi_j,\varphi_k)}{\partial_x u_0}
= \tilde{p}_j\,\bar{a}_j + \tilde{Q}_{jk}\,\bar{a}_j\bar{a}_k ,
\]
the coefficients being exactly those of~\eqref{eq:app-sr-pq} ($\tilde Q$ symmetric because $B$ is).
\end{proof}

\begin{corollary}[Exact recovery of the symmetry-reduced Galerkin model]
\label{cor:app-sr}
Assume \textnormal{(A1)--(A3)} and that the re-projected states $\{\bar{a}(t_m)\}$ make the shared regression matrix with rows $[\,\bar{a}_j(t_m),\ \bar{a}_j(t_m)\bar{a}_k(t_m)\,]$ of full column rank, which requires
\begin{equation}
\label{eq:app-sr-rank}
K \;\ge\; n + \frac{n(n+1)}{2}
\end{equation}
re-projected samples (duplicate quadratic monomials identified).  Then both re-projected least-squares problems---\eqref{eq:app-sr-ls-f} with target $\bar{v}(t_m)$ and~\eqref{eq:app-sr-ls-c} with target $\bar{g}(t_m)$---have unique minimizers attaining zero residual, equal to the symmetry-reduced Galerkin coefficients, $\hat{A}=\tilde{A}$, $\hat{B}=\tilde{B}$, $\hat{p}=\tilde{p}$, $\hat{Q}=\tilde{Q}$.  Together with the basis quantities $\tilde{C}_{ij}$ and $s_j$, the assembled symmetry-reduced operator inference ROM then coincides exactly with the symmetry-reduced Galerkin ROM.
\end{corollary}

\begin{proof}
Both regressions share the feature map $a\mapsto[\,a_j;\,a_j a_k\,]$, so the rank requirement is the single condition~\eqref{eq:app-sr-rank}.  By Proposition~\ref{prop:app-sr} the symmetry-reduced Galerkin coefficients drive each nonnegative loss to zero and are global minimizers; full column rank makes each unique.  The remaining coefficients $\tilde{C}_{ij}$ and $s_j$ are basis quantities, identical in both models.
\end{proof}

\begin{remark}[What re-projection does and does not touch]
\label{rem:app-sr}
Only the two $f$-dependent targets are re-projected; the shift-coupling coefficients $\tilde{C}_{ij}$ and the denominator coefficients $s_j$ are known from the basis and template and require neither learning nor re-projection.  A single re-projected dataset $\{\bar{a}(t_m)\}$---one query of $f$ at each lifted state $P u_{\mathrm{fit}}(t_m)$---simultaneously supplies the targets $\bar{v}(t_m)$ and $\bar{g}(t_m)$ for both regressions, so the entire symmetry-reduced Galerkin ROM is recovered from the same data under the single condition~\eqref{eq:app-sr-rank}.
\end{remark}

\end{appendices}



\begin{thebibliography}{49}
\ifx \bisbn   \undefined \def \bisbn  #1{ISBN #1}\fi
\ifx \binits  \undefined \def \binits#1{#1}\fi
\ifx \bauthor  \undefined \def \bauthor#1{#1}\fi
\ifx \batitle  \undefined \def \batitle#1{#1}\fi
\ifx \bjtitle  \undefined \def \bjtitle#1{#1}\fi
\ifx \bvolume  \undefined \def \bvolume#1{\textbf{#1}}\fi
\ifx \byear  \undefined \def \byear#1{#1}\fi
\ifx \bissue  \undefined \def \bissue#1{#1}\fi
\ifx \bfpage  \undefined \def \bfpage#1{#1}\fi
\ifx \blpage  \undefined \def \blpage #1{#1}\fi
\ifx \burl  \undefined \def \burl#1{\textsf{#1}}\fi
\ifx \doiurl  \undefined \def \doiurl#1{\url{https://doi.org/#1}}\fi
\ifx \betal  \undefined \def \betal{\textit{et al.}}\fi
\ifx \binstitute  \undefined \def \binstitute#1{#1}\fi
\ifx \binstitutionaled  \undefined \def \binstitutionaled#1{#1}\fi
\ifx \bctitle  \undefined \def \bctitle#1{#1}\fi
\ifx \beditor  \undefined \def \beditor#1{#1}\fi
\ifx \bpublisher  \undefined \def \bpublisher#1{#1}\fi
\ifx \bbtitle  \undefined \def \bbtitle#1{#1}\fi
\ifx \bedition  \undefined \def \bedition#1{#1}\fi
\ifx \bseriesno  \undefined \def \bseriesno#1{#1}\fi
\ifx \blocation  \undefined \def \blocation#1{#1}\fi
\ifx \bsertitle  \undefined \def \bsertitle#1{#1}\fi
\ifx \bsnm \undefined \def \bsnm#1{#1}\fi
\ifx \bsuffix \undefined \def \bsuffix#1{#1}\fi
\ifx \bparticle \undefined \def \bparticle#1{#1}\fi
\ifx \barticle \undefined \def \barticle#1{#1}\fi
\bibcommenthead
\ifx \bconfdate \undefined \def \bconfdate #1{#1}\fi
\ifx \botherref \undefined \def \botherref #1{#1}\fi
\ifx \url \undefined \def \url#1{\textsf{#1}}\fi
\ifx \bchapter \undefined \def \bchapter#1{#1}\fi
\ifx \bbook \undefined \def \bbook#1{#1}\fi
\ifx \bcomment \undefined \def \bcomment#1{#1}\fi
\ifx \oauthor \undefined \def \oauthor#1{#1}\fi
\ifx \citeauthoryear \undefined \def \citeauthoryear#1{#1}\fi
\ifx \endbibitem  \undefined \def \endbibitem {}\fi
\ifx \bconflocation  \undefined \def \bconflocation#1{#1}\fi
\ifx \arxivurl  \undefined \def \arxivurl#1{\textsf{#1}}\fi
\csname PreBibitemsHook\endcsname

\bibitem[\protect\citeauthoryear{Wedin and Kerswell}{2004}]{Wedin2004}
\begin{barticle}
\bauthor{\bsnm{Wedin}, \binits{H.}},
\bauthor{\bsnm{Kerswell}, \binits{R.R.}}:
\batitle{Exact coherent structures in pipe flow: travelling wave solutions}.
\bjtitle{J. Fluid Mech.}
\bvolume{508},
\bfpage{333}--\blpage{371}
(\byear{2004})
\doiurl{10.1017/S0022112004009346}
\end{barticle}
\endbibitem

\bibitem[\protect\citeauthoryear{Avila et~al.}{2013}]{Avila2013}
\begin{barticle}
\bauthor{\bsnm{Avila}, \binits{M.}},
\bauthor{\bsnm{Mellibovsky}, \binits{F.}},
\bauthor{\bsnm{Roland}, \binits{N.}},
\bauthor{\bsnm{Hof}, \binits{B.}}:
\batitle{Streamwise-localized solutions at the onset of turbulence in pipe flow}.
\bjtitle{Phys. Rev. Lett.}
\bvolume{110},
\bfpage{224502}
(\byear{2013})
\doiurl{10.1103/PhysRevLett.110.224502}
\end{barticle}
\endbibitem

\bibitem[\protect\citeauthoryear{Willis et~al.}{2013}]{Willis2013}
\begin{barticle}
\bauthor{\bsnm{Willis}, \binits{A.P.}},
\bauthor{\bsnm{Cvitanović}, \binits{P.}},
\bauthor{\bsnm{Avila}, \binits{M.}}:
\batitle{Revealing the state space of turbulent pipe flow by symmetry reduction}.
\bjtitle{J. Fluid Mech.}
\bvolume{721},
\bfpage{514}--\blpage{540}
(\byear{2013})
\doiurl{10.1017/jfm.2013.75}
\end{barticle}
\endbibitem

\bibitem[\protect\citeauthoryear{Willis et~al.}{2016}]{Willis2016}
\begin{barticle}
\bauthor{\bsnm{Willis}, \binits{A.P.}},
\bauthor{\bsnm{Short}, \binits{K.Y.}},
\bauthor{\bsnm{Cvitanović}, \binits{P.}}:
\batitle{Symmetry reduction in high dimensions, illustrated in a turbulent pipe}.
\bjtitle{Phys. Rev. E}
\bvolume{93},
\bfpage{022204}
(\byear{2016})
\doiurl{10.1103/PhysRevE.93.022204}
\end{barticle}
\endbibitem

\bibitem[\protect\citeauthoryear{Ritter et~al.}{2016}]{Ritter2016}
\begin{barticle}
\bauthor{\bsnm{Ritter}, \binits{P.}},
\bauthor{\bsnm{Mellibovsky}, \binits{F.}},
\bauthor{\bsnm{Avila}, \binits{M.}}:
\batitle{Emergence of spatio-temporal dynamics from exact coherent solutions in pipe flow}.
\bjtitle{New J. Phys.}
\bvolume{18}(\bissue{8}),
\bfpage{083031}
(\byear{2016})
\doiurl{10.1088/1367-2630/18/8/083031}
\end{barticle}
\endbibitem

\bibitem[\protect\citeauthoryear{Budanur et~al.}{2017}]{Budanur2017}
\begin{barticle}
\bauthor{\bsnm{Budanur}, \binits{N.B.}},
\bauthor{\bsnm{Short}, \binits{K.Y.}},
\bauthor{\bsnm{Farazmand}, \binits{M.}},
\bauthor{\bsnm{Willis}, \binits{A.P.}},
\bauthor{\bsnm{Cvitanović}, \binits{P.}}:
\batitle{Relative periodic orbits form the backbone of turbulent pipe flow}.
\bjtitle{J. Fluid Mech.}
\bvolume{833},
\bfpage{274}--\blpage{301}
(\byear{2017})
\doiurl{10.1017/jfm.2017.699}
\end{barticle}
\endbibitem

\bibitem[\protect\citeauthoryear{Sirovich}{1987}]{Sirovich1987}
\begin{barticle}
\bauthor{\bsnm{Sirovich}, \binits{L.}}:
\batitle{Turbulence and the dynamics of coherent structures part i: coherent structures}.
\bjtitle{Q. Appl. Math.}
\bvolume{45}(\bissue{3}),
\bfpage{561}--\blpage{571}
(\byear{1987})
{\href{https://arxiv.org/abs/http://www.jstor.org/stable/43637457}{{http://www.jstor.org/stable/43637457}}}
\end{barticle}
\endbibitem

\bibitem[\protect\citeauthoryear{Berkooz et~al.}{1993}]{Berkooz1993}
\begin{barticle}
\bauthor{\bsnm{Berkooz}, \binits{G.}},
\bauthor{\bsnm{Holmes}, \binits{P.}},
\bauthor{\bsnm{Lumley}, \binits{J.L.}}:
\batitle{The proper orthogonal decomposition in the analysis of turbulent flows}.
\bjtitle{Annu. Rev. Fluid Mech.}
\bvolume{25},
\bfpage{539}--\blpage{575}
(\byear{1993})
\doiurl{10.1146/annurev.fl.25.010193.002543}
\end{barticle}
\endbibitem

\bibitem[\protect\citeauthoryear{Philip et~al.}{2012}]{Holmes2012}
\begin{bbook}
\bauthor{\bsnm{Philip}, \binits{H.}},
\bauthor{\bsnm{Lumley}, \binits{J.L.}},
\bauthor{\bsnm{Berkooz}, \binits{G.}},
\bauthor{\bsnm{Rowley}, \binits{C.W.}}:
\bbtitle{Turbulence, {C}oherent {S}tructures, {D}ynamical {S}ystems And {S}ymmetry},
\bedition{2}nd edn.
\bpublisher{Cambridge University Press},
\blocation{Cambridge}
(\byear{2012}).
\doiurl{10.1017/CBO9780511622700}
\end{bbook}
\endbibitem

\bibitem[\protect\citeauthoryear{Greif and Urban}{2019}]{Greif2019}
\begin{barticle}
\bauthor{\bsnm{Greif}, \binits{C.}},
\bauthor{\bsnm{Urban}, \binits{K.}}:
\batitle{Decay of the {K}olmogorov {N}-width for wave problems}.
\bjtitle{Appl. Math. Lett.}
\bvolume{96},
\bfpage{216}--\blpage{222}
(\byear{2019})
\doiurl{10.1016/j.aml.2019.05.013}
\end{barticle}
\endbibitem

\bibitem[\protect\citeauthoryear{Kirby and Armbruster}{1992}]{Kirby1992}
\begin{barticle}
\bauthor{\bsnm{Kirby}, \binits{M.}},
\bauthor{\bsnm{Armbruster}, \binits{D.}}:
\batitle{Reconstructing phase space from {PDE} simulations}.
\bjtitle{Z. angew. Math. Phys.}
\bvolume{43},
\bfpage{999}--\blpage{1022}
(\byear{1992})
\doiurl{10.1007/BF00916425}
\end{barticle}
\endbibitem

\bibitem[\protect\citeauthoryear{Rowley and Marsden}{2000}]{Rowley2000}
\begin{barticle}
\bauthor{\bsnm{Rowley}, \binits{C.W.}},
\bauthor{\bsnm{Marsden}, \binits{J.E.}}:
\batitle{Reconstruction equations and the {K}arhunen-lo{\`e}ve expansion for systems with symmetry}.
\bjtitle{Physica D}
\bvolume{142},
\bfpage{1}--\blpage{19}
(\byear{2000})
\doiurl{10.1016/S0167-2789(00)00042-7}
\end{barticle}
\endbibitem

\bibitem[\protect\citeauthoryear{Rowley et~al.}{2003}]{Rowley-nl03}
\begin{barticle}
\bauthor{\bsnm{Rowley}, \binits{C.W.}},
\bauthor{\bsnm{Kevrekidis}, \binits{I.G.}},
\bauthor{\bsnm{Marsden}, \binits{J.E.}},
\bauthor{\bsnm{Lust}, \binits{K.}}:
\batitle{Reduction and reconstruction for self-similar dynamical systems}.
\bjtitle{Nonlinearity}
\bvolume{16}(\bissue{4}),
\bfpage{1257}
(\byear{2003})
\doiurl{10.1088/0951-7715/16/4/304}
\end{barticle}
\endbibitem

\bibitem[\protect\citeauthoryear{Budanur et~al.}{2015}]{Budanur2015}
\begin{barticle}
\bauthor{\bsnm{Budanur}, \binits{N.B.}},
\bauthor{\bsnm{Cvitanović}, \binits{P.}},
\bauthor{\bsnm{Davidchack}, \binits{R.L.}},
\bauthor{\bsnm{Siminos}, \binits{E.}}:
\batitle{Reduction of {SO}(2) symmetry for spatially extended dynamical systems}.
\bjtitle{Phys. Rev. Lett.}
\bvolume{114},
\bfpage{084102}
(\byear{2015})
\doiurl{10.1103/PhysRevLett.114.084102}
\end{barticle}
\endbibitem

\bibitem[\protect\citeauthoryear{Marsden et~al.}{1990}]{Marsden1990}
\begin{botherref}
\oauthor{\bsnm{Marsden}, \binits{J.E.}},
\oauthor{\bsnm{Montgomery}, \binits{R.}},
\oauthor{\bsnm{Ratiu}, \binits{T.S.}}:
Reduction, symmetry and phases in mechanics.
Mem. Am. Math. Soc.
\textbf{88}(436)
(1990)
\doiurl{10.1090/memo/0436}
\end{botherref}
\endbibitem

\bibitem[\protect\citeauthoryear{Marsden et~al.}{2000}]{Marsden2000}
\begin{barticle}
\bauthor{\bsnm{Marsden}, \binits{J.E.}},
\bauthor{\bsnm{Ratiu}, \binits{T.S.}},
\bauthor{\bsnm{Scheurle}, \binits{J.}}:
\batitle{Reduction theory and the {Lagrange}–{Routh} equations}.
\bjtitle{J. Math. Phys.}
\bvolume{41}(\bissue{6}),
\bfpage{3379}--\blpage{3429}
(\byear{2000})
\doiurl{10.1063/1.533317}
\end{barticle}
\endbibitem

\bibitem[\protect\citeauthoryear{Beyn and Th{\"u}mmler}{2004}]{Beyn2004}
\begin{barticle}
\bauthor{\bsnm{Beyn}, \binits{W.}},
\bauthor{\bsnm{Th{\"u}mmler}, \binits{V.}}:
\batitle{Freezing solutions of equivariant evolution equations}.
\bjtitle{SIAM J. Appl. Dyn. Syst.}
\bvolume{3}(\bissue{2}),
\bfpage{85}--\blpage{116}
(\byear{2004})
\doiurl{10.1137/030600515}
\end{barticle}
\endbibitem

\bibitem[\protect\citeauthoryear{Black et~al.}{2020}]{Black2020}
\begin{barticle}
\bauthor{\bsnm{Black}, \binits{F.}},
\bauthor{\bsnm{Schulze}, \binits{P.}},
\bauthor{\bsnm{Unger}, \binits{B.}}:
\batitle{Projection-based model reduction with dynamically transformed modes}.
\bjtitle{ESAIM: Math. Model. Numer. Anal.}
\bvolume{54}(\bissue{6}),
\bfpage{2011}--\blpage{2043}
(\byear{2020})
\doiurl{10.1051/m2an/2020046}
\end{barticle}
\endbibitem

\bibitem[\protect\citeauthoryear{Reiss et~al.}{2018}]{Reiss2019}
\begin{barticle}
\bauthor{\bsnm{Reiss}, \binits{J.}},
\bauthor{\bsnm{Schulze}, \binits{P.}},
\bauthor{\bsnm{Sesterhenn}, \binits{J.}},
\bauthor{\bsnm{Mehrmann}, \binits{V.}}:
\batitle{The shifted proper orthogonal decomposition: A mode decomposition for multiple transport phenomena}.
\bjtitle{SIAM J. Sci. Comput.}
\bvolume{40}(\bissue{3}),
\bfpage{1322}--\blpage{1344}
(\byear{2018})
\doiurl{10.1137/17M1140571}
\end{barticle}
\endbibitem

\bibitem[\protect\citeauthoryear{Mendible et~al.}{2020}]{Mendible2020}
\begin{barticle}
\bauthor{\bsnm{Mendible}, \binits{A.}},
\bauthor{\bsnm{Brunton}, \binits{S.L.}},
\bauthor{\bsnm{Aravkin}, \binits{A.Y.}},
\bauthor{\bsnm{Lowrie}, \binits{W.}},
\bauthor{\bsnm{Nathan~Kutz}, \binits{J.}}:
\batitle{Dimensionality reduction and reduced-order modeling for traveling wave physics}.
\bjtitle{Theor. Comput. Fluid Dyn.}
\bvolume{34},
\bfpage{385}--\blpage{400}
(\byear{2020})
\doiurl{10.1007/s00162-020-00529-9}
\end{barticle}
\endbibitem

\bibitem[\protect\citeauthoryear{Lee and Carlberg}{2020}]{LeeCarlberg2020}
\begin{barticle}
\bauthor{\bsnm{Lee}, \binits{K.}},
\bauthor{\bsnm{Carlberg}, \binits{K.T.}}:
\batitle{Model reduction of dynamical systems on nonlinear manifolds using deep convolutional autoencoders}.
\bjtitle{J. Comput. Phys}
\bvolume{404},
\bfpage{108973}
(\byear{2020})
\doiurl{10.1016/j.jcp.2019.108973}
\end{barticle}
\endbibitem

\bibitem[\protect\citeauthoryear{Peherstorfer}{2020}]{Peherstorfer2020_online_adaptive_basis}
\begin{barticle}
\bauthor{\bsnm{Peherstorfer}, \binits{B.}}:
\batitle{Model reduction for transport-dominated problems via online adaptive bases and adaptive sampling}.
\bjtitle{SIAM J. Sci. Comput.}
\bvolume{42}(\bissue{5}),
\bfpage{2803}--\blpage{2836}
(\byear{2020})
\doiurl{10.1137/19M1257275}
\end{barticle}
\endbibitem

\bibitem[\protect\citeauthoryear{Rowley et~al.}{2009}]{Rowley2009}
\begin{barticle}
\bauthor{\bsnm{Rowley}, \binits{C.W.}},
\bauthor{\bsnm{Mezi\'c}, \binits{I.}},
\bauthor{\bsnm{Bagheri}, \binits{S.}},
\bauthor{\bsnm{Schlatter}, \binits{P.}},
\bauthor{\bsnm{Henningson}, \binits{D.S.}}:
\batitle{Spectral analysis of nonlinear flows}.
\bjtitle{J. Fluid. Mech.}
\bvolume{641},
\bfpage{115}--\blpage{127}
(\byear{2009})
\doiurl{10.1017/S0022112009992059}
\end{barticle}
\endbibitem

\bibitem[\protect\citeauthoryear{Schmid}{2010}]{Schmid2010}
\begin{barticle}
\bauthor{\bsnm{Schmid}, \binits{P.J.}}:
\batitle{Dynamic mode decomposition of numerical and experimental data}.
\bjtitle{J. Fluid Mech.}
\bvolume{656},
\bfpage{5}--\blpage{28}
(\byear{2010})
\doiurl{10.1017/S0022112010001217}
\end{barticle}
\endbibitem

\bibitem[\protect\citeauthoryear{Tu et~al.}{2014}]{Tu2012}
\begin{barticle}
\bauthor{\bsnm{Tu}, \binits{J.H.}},
\bauthor{\bsnm{Rowley}, \binits{C.W.}},
\bauthor{\bsnm{Luchtenburg}, \binits{D.M.}},
\bauthor{\bsnm{Brunton}, \binits{S.L.}},
\bauthor{\bsnm{Nathan~Kutz}, \binits{J.}}:
\batitle{On dynamic mode decomposition: Theory and applications}.
\bjtitle{J. Comput. Dyn.}
\bvolume{1}(\bissue{2}),
\bfpage{391}--\blpage{421}
(\byear{2014})
\doiurl{10.3934/jcd.2014.1.391}
\end{barticle}
\endbibitem

\bibitem[\protect\citeauthoryear{Peherstorfer and Willcox}{2016}]{Peherstorfer2016}
\begin{barticle}
\bauthor{\bsnm{Peherstorfer}, \binits{B.}},
\bauthor{\bsnm{Willcox}, \binits{K.}}:
\batitle{Data-driven operator inference for nonintrusive projection-based model reduction}.
\bjtitle{Comput. Methods Appl. Mech. Eng.}
\bvolume{306},
\bfpage{196}--\blpage{215}
(\byear{2016})
\doiurl{10.1016/j.cma.2016.03.025}
\end{barticle}
\endbibitem

\bibitem[\protect\citeauthoryear{Qian et~al.}{2020}]{Qian2020}
\begin{barticle}
\bauthor{\bsnm{Qian}, \binits{E.}},
\bauthor{\bsnm{Kramer}, \binits{B.}},
\bauthor{\bsnm{Peherstorfer}, \binits{B.}},
\bauthor{\bsnm{Willcox}, \binits{K.}}:
\batitle{{L}ift \& {L}earn: {P}hysics-informed machine learning for large-scale nonlinear dynamical systems}.
\bjtitle{Physica D}
\bvolume{406},
\bfpage{132401}
(\byear{2020})
\doiurl{10.1016/j.physd.2020.132401}
\end{barticle}
\endbibitem

\bibitem[\protect\citeauthoryear{Kramer et~al.}{2024}]{Kramer2024}
\begin{barticle}
\bauthor{\bsnm{Kramer}, \binits{B.}},
\bauthor{\bsnm{Peherstorfer}, \binits{B.}},
\bauthor{\bsnm{Willcox}, \binits{K.E.}}:
\batitle{Learning nonlinear reduced models from data with operator inference}.
\bjtitle{Annu. Rev. Fluid Mech.}
\bvolume{56},
\bfpage{521}--\blpage{548}
(\byear{2024})
\doiurl{10.1146/annurev-fluid-121021-025220}
\end{barticle}
\endbibitem

\bibitem[\protect\citeauthoryear{Baddoo et~al.}{2023}]{Baddoo2023}
\begin{barticle}
\bauthor{\bsnm{Baddoo}, \binits{P.J.}},
\bauthor{\bsnm{Herrmann}, \binits{B.}},
\bauthor{\bsnm{McKeon}, \binits{B.J.}},
\bauthor{\bsnm{Nathan~Kutz}, \binits{J.}},
\bauthor{\bsnm{Brunton}, \binits{S.L.}}:
\batitle{Physics-informed dynamic mode decomposition}.
\bjtitle{Proc. R. Soc. A}
\bvolume{479}(\bissue{2271}),
\bfpage{20220576}
(\byear{2023})
\doiurl{10.1098/rspa.2022.0576}
\end{barticle}
\endbibitem

\bibitem[\protect\citeauthoryear{Sesterhenn and Shahirpour}{2019}]{Sesterhenn2019}
\begin{barticle}
\bauthor{\bsnm{Sesterhenn}, \binits{J.}},
\bauthor{\bsnm{Shahirpour}, \binits{A.}}:
\batitle{A characteristic dynamic mode decomposition}.
\bjtitle{Theor. Comput. Fluid Dyn.}
\bvolume{33},
\bfpage{281}--\blpage{305}
(\byear{2019})
\doiurl{10.1007/s00162-019-00494-y}
\end{barticle}
\endbibitem

\bibitem[\protect\citeauthoryear{Marensi et~al.}{2023}]{Marensi2023}
\begin{barticle}
\bauthor{\bsnm{Marensi}, \binits{E.}},
\bauthor{\bsnm{Yalnız}, \binits{G.}},
\bauthor{\bsnm{Hof}, \binits{B.}},
\bauthor{\bsnm{Budanur}, \binits{N.B.}}:
\batitle{Symmetry-reduced dynamic mode decomposition of near-wall turbulence}.
\bjtitle{J. Fluid Mech.}
\bvolume{954},
\bfpage{10}
(\byear{2023})
\doiurl{10.1017/jfm.2022.1001}
\end{barticle}
\endbibitem

\bibitem[\protect\citeauthoryear{Engel et~al.}{2025}]{Engel2025}
\begin{barticle}
\bauthor{\bsnm{Engel}, \binits{M.}},
\bauthor{\bsnm{Ashtari}, \binits{O.}},
\bauthor{\bsnm{Schneider}, \binits{T.M.}},
\bauthor{\bsnm{Linkmann}, \binits{M.}}:
\batitle{Search for unstable relative periodic orbits in channel flow using symmetry-reduced dynamic mode decomposition}.
\bjtitle{J. Fluid. Mech.}
\bvolume{1013},
\bfpage{45}
(\byear{2025})
\doiurl{10.1017/jfm.2025.10255}
\end{barticle}
\endbibitem

\bibitem[\protect\citeauthoryear{Issan and Kramer}{2023}]{Issan2023}
\begin{barticle}
\bauthor{\bsnm{Issan}, \binits{O.}},
\bauthor{\bsnm{Kramer}, \binits{B.}}:
\batitle{Predicting solar wind streams from the inner-heliosphere to {E}arth via shifted operator inference}.
\bjtitle{J. Comput. Phys.}
\bvolume{473},
\bfpage{111689}
(\byear{2023})
\doiurl{10.1016/j.jcp.2022.111689}
\end{barticle}
\endbibitem

\bibitem[\protect\citeauthoryear{Burela et~al.}{2025}]{Burela2025}
\begin{barticle}
\bauthor{\bsnm{Burela}, \binits{S.}},
\bauthor{\bsnm{Krah}, \binits{P.}},
\bauthor{\bsnm{Reiss}, \binits{J.}}:
\batitle{Parametric model order reduction for a wildland fire model via the shifted {POD}-based deep learning method}.
\bjtitle{Adv. Comput. Math.}
\bvolume{51}(\bissue{1}),
\bfpage{9}
(\byear{2025})
\doiurl{10.1007/s10444-025-10220-4}
\end{barticle}
\endbibitem

\bibitem[\protect\citeauthoryear{Cagniart et~al.}{2019}]{Cagniart2019}
\begin{bbook}
\bauthor{\bsnm{Cagniart}, \binits{N.}},
\bauthor{\bsnm{Maday}, \binits{Y.}},
\bauthor{\bsnm{Stamm}, \binits{B.}}:
In: \beditor{\bsnm{Chetverushkin}, \binits{B.N.}},
\beditor{\bsnm{Fitzgibbon}, \binits{W.}},
\beditor{\bsnm{Kuznetsov}, \binits{Y.A.}},
\beditor{\bsnm{Neittaanmäki}, \binits{P.}},
\beditor{\bsnm{Periaux}, \binits{J.}},
\beditor{\bsnm{Pironneau}, \binits{O.}} (eds.)
\bbtitle{Model order reduction for problems with large convection effects}.
\bsertitle{Computational Methods in Applied Sciences},
vol. \bseriesno{47},
pp. \bfpage{131}--\blpage{150}.
\bpublisher{Springer},
\blocation{Cham, Switzerland}
(\byear{2019}).
\doiurl{10.1007/978-3-319-78325-3_10}
\end{bbook}
\endbibitem

\bibitem[\protect\citeauthoryear{Peherstorfer}{2020}]{Peherstorfer2020}
\begin{barticle}
\bauthor{\bsnm{Peherstorfer}, \binits{B.}}:
\batitle{Sampling low-dimensional markovian dynamics for preasymptotically recovering reduced models from data with operator inference}.
\bjtitle{SIAM J. Sci. Comput.}
\bvolume{42}(\bissue{5}),
\bfpage{3489}--\blpage{3515}
(\byear{2020})
\doiurl{10.1137/19M1292448}
\end{barticle}
\endbibitem

\bibitem[\protect\citeauthoryear{Kevrekidis et~al.}{1990}]{Kevrekidis1990}
\begin{barticle}
\bauthor{\bsnm{Kevrekidis}, \binits{I.G.}},
\bauthor{\bsnm{Nicolaenko}, \binits{B.}},
\bauthor{\bsnm{Scovel}, \binits{J.C.}}:
\batitle{Back in the saddle again: A computer assisted study of the {K}uramoto–{S}ivashinsky equation}.
\bjtitle{SIAM J. Appl. Math.}
\bvolume{50}(\bissue{3}),
\bfpage{760}--\blpage{790}
(\byear{1990})
\doiurl{10.1137/0150045}
\end{barticle}
\endbibitem

\bibitem[\protect\citeauthoryear{Orszag}{1971}]{Orszag1971}
\begin{barticle}
\bauthor{\bsnm{Orszag}, \binits{S.A.}}:
\batitle{On the elimination of aliasing in finite-difference schemes by filtering high-wavenumber components}.
\bjtitle{J. Atmos. Sci.}
\bvolume{28}(\bissue{6}),
\bfpage{1074}
(\byear{1971})
\end{barticle}
\endbibitem

\bibitem[\protect\citeauthoryear{Zang and Hussaini}{}]{Zhang1985}
\begin{botherref}
\oauthor{\bsnm{Zang}, \binits{T.}},
\oauthor{\bsnm{Hussaini}, \binits{M.}}:
Numerical experiments on subcritical transition mechanisms.
AIAA 23rd {A}erospace {S}ciences {M}eeting,
1985--0296
\doiurl{10.2514/6.1985-296}
\end{botherref}
\endbibitem

\bibitem[\protect\citeauthoryear{Peyret}{2002}]{Peyret2002}
\begin{bbook}
\bauthor{\bsnm{Peyret}, \binits{R.}}:
\bbtitle{Spectral {M}ethods for {I}ncompressible {V}iscous {F}low}.
\bpublisher{Springer},
\blocation{New York}
(\byear{2002}).
\doiurl{10.1007/978-1-4757-6557-1}
\end{bbook}
\endbibitem

\bibitem[\protect\citeauthoryear{Hairer and Wanner}{1996}]{Hairer1996}
\begin{bbook}
\bauthor{\bsnm{Hairer}, \binits{E.}},
\bauthor{\bsnm{Wanner}, \binits{G.}}:
\bbtitle{Solving Ordinary Differential Equations {II}: Stiff and Differential-Algebraic Problems},
\bedition{2nd} edn.
\bsertitle{Springer Series in Computational Mathematics},
vol. \bseriesno{14}.
\bpublisher{Springer},
\blocation{Berlin, Heidelberg}
(\byear{1996}).
\doiurl{10.1007/978-3-642-05221-7}
\end{bbook}
\endbibitem

\bibitem[\protect\citeauthoryear{Swischuk et~al.}{2020}]{Swischuk2020}
\begin{barticle}
\bauthor{\bsnm{Swischuk}, \binits{R.}},
\bauthor{\bsnm{Kramer}, \binits{B.}},
\bauthor{\bsnm{Huang}, \binits{C.}},
\bauthor{\bsnm{Willcox}, \binits{K.}}:
\batitle{Learning physics-based reduced-order models for a single-injector combustion process}.
\bjtitle{AIAA J.}
\bvolume{58}(\bissue{6}),
\bfpage{2658}--\blpage{2672}
(\byear{2020})
\doiurl{10.2514/1.J058943}
\end{barticle}
\endbibitem

\bibitem[\protect\citeauthoryear{Jain et~al.}{}]{Jain2021}
\begin{botherref}
\oauthor{\bsnm{Jain}, \binits{P.}},
\oauthor{\bsnm{S.}, \binits{M.}},
\oauthor{\bsnm{Kramer}, \binits{B.}}:
Performance comparison of data-driven reduced models for a single-injector combustion process.
\doiurl{10.2514/6.2021-3633}
\end{botherref}
\endbibitem

\bibitem[\protect\citeauthoryear{McQuarrie et~al.}{2021}]{McQuarrie2021}
\begin{botherref}
\oauthor{\bsnm{McQuarrie}, \binits{S.A.}},
\oauthor{\bsnm{H.}, \binits{C.}},
\oauthor{\bsnm{Willcox}, \binits{K.E.}}:
Data-driven reduced-order models via regularised operator inference for a single-injector combustion process.
J. R. Soc. N. Z.
\textbf{51}(2)
(2021)
\doiurl{10.1080/03036758.2020.1863237}
\end{botherref}
\endbibitem

\bibitem[\protect\citeauthoryear{Sawant et~al.}{2023}]{Sawant2021}
\begin{barticle}
\bauthor{\bsnm{Sawant}, \binits{N.}},
\bauthor{\bsnm{Kramer}, \binits{B.}},
\bauthor{\bsnm{Peherstorfer}, \binits{B.}}:
\batitle{Physics-informed regularization and structure preservation for learning stable reduced models from data with operator inference}.
\bjtitle{Comput. Methods Appl. Mech. Eng.}
\bvolume{404},
\bfpage{115836}
(\byear{2023})
\doiurl{10.1016/j.cma.2022.115836}
\end{barticle}
\endbibitem

\bibitem[\protect\citeauthoryear{Farcas et~al.}{}]{Farcas2023}
\begin{botherref}
\oauthor{\bsnm{Farcas}, \binits{I.}},
\oauthor{\bsnm{R.}, \binits{G.}},
\oauthor{\bsnm{Munipalli}, \binits{R.}},
\oauthor{\bsnm{Willcox}, \binits{K.E.}}:
Parametric non-intrusive reduced-order models via operator inference for large-scale rotating detonation engine simulations.
\doiurl{10.2514/6.2023-0172}
\end{botherref}
\endbibitem

\bibitem[\protect\citeauthoryear{Padovan et~al.}{2024}]{Padovan2024}
\begin{barticle}
\bauthor{\bsnm{Padovan}, \binits{A.}},
\bauthor{\bsnm{Vollmer}, \binits{B.}},
\bauthor{\bsnm{Bodony}, \binits{D.J.}}:
\batitle{Data-driven model reduction via non-intrusive optimization of projection operators and reduced-order dynamics}.
\bjtitle{SIAM J. Appl. Dyn. Syst.}
\bvolume{23}(\bissue{4}),
\bfpage{3052}--\blpage{3076}
(\byear{2024})
\doiurl{10.1137/24M1628414}
\end{barticle}
\endbibitem

\bibitem[\protect\citeauthoryear{McQuarrie et~al.}{2023}]{McQuarrie2023}
\begin{barticle}
\bauthor{\bsnm{McQuarrie}, \binits{S.A.}},
\bauthor{\bsnm{Khodabakhshi}, \binits{P.}},
\bauthor{\bsnm{Willcox}, \binits{K.E.}}:
\batitle{Nonintrusive reduced-order models for parametric partial differential equations via data-driven operator inference}.
\bjtitle{SIAM J. Sci. Comput.}
\bvolume{45}(\bissue{4}),
\bfpage{1917}--\blpage{1946}
(\byear{2023})
\doiurl{10.1137/21M1452810}
\end{barticle}
\endbibitem

\bibitem[\protect\citeauthoryear{Kramer et~al.}{2017}]{Kramer2017}
\begin{barticle}
\bauthor{\bsnm{Kramer}, \binits{B.}},
\bauthor{\bsnm{Peherstorfer}, \binits{B.}},
\bauthor{\bsnm{Willcox}, \binits{K.}}:
\batitle{Feedback control for systems with uncertain parameters using online-adaptive reduced models}.
\bjtitle{SIAM J. Appl. Dyn. Syst.}
\bvolume{16}(\bissue{3}),
\bfpage{1563}--\blpage{1586}
(\byear{2017})
\doiurl{10.1137/16M1088958}
\end{barticle}
\endbibitem

\end{thebibliography}

\end{document}